\documentclass[a4paper,12pt]{article}
\usepackage{amsmath,epsfig,pstricks,pst-node}
\oddsidemargin -0.15in\newtheorem{Theorem}{Theorem}[section]
\newtheorem{Corollary}[Theorem]{Corollary}

\newtheorem{Definition}[Theorem]{Definition}

\newcommand{\beq}{\begin{equation}}
\newcommand{\eeq}{\end{equation}}

\newenvironment{fullfigure}[2]
    {\begin{figure}[htb]\begin{center}\def\ffa{#1}\def\ffb{#2}}
    {\vspace{\baselineskip}\caption{\ffb.}\label{\ffa}\end{center}\end{figure}}

\def\ssc{\scriptscriptstyle}

\def\ov{\overline}

\def\ph1{{\phantom1}}
\def\ph-{{\phantom-}}

\def\con{{\rm con}}
\def\col{{\rm col}}

\def\wgt{{\rm wgt}}
\def\var{{\rm var}}
\def\hgt{{\rm hgt}}
\def\str{{\rm str}}

\def\bar{{\rm bar}}

\def\neg{{\rm neg}}
\def\row{{\rm row}}
\def\ssi{{\rm ssi}}
\def\we{{\rm we}}
\def\ns{{\rm ns}}
\def\ne{{\rm ne}}
\def\sw{{\rm sw}}
\def\nw{{\rm nw}}
\def\se{{\rm se}}

\def\={\! = \!}
\def\p{\! + \!}
\def\m{\! - \!}

\def\0{\phantom{0}}
\def\o{\circ}

\def\emdash{\hbox{---}}
\def\phW{\phantom{\,\emdash H}}
\def\phE{\phantom{H\emdash\,}}

\def\la{\leftarrow}
\def\ra{\rightarrow}
\def\ua{\uparrow}
\def\da{\downarrow}
\def\HOH{H\ra O\la H}
\def\HO{H\ra O}
\def\OH{O\la H}

\def\mystrut{\hbox{\vrule height10.6pt depth2pt width0pt}}
\def\mybox{\hbox to 8.9pt}
\def\norulefill{\leaders\hrule height0pt\hfill}
\def\nr#1{\multispan{#1}\norulefill}
\def\hr#1{\multispan{#1}\hrulefill}

\begin{document}
\title{U--turn alternating sign matrices, symplectic shifted tableaux
and their weighted enumeration}
\author{A M Hamel\thanks{e-mail: ahamel@wlu.ca} \\
Department of Physics and Computer Science,\\
Wilfrid Laurier University, 
Waterloo, Ontario N2L 3C5, Canada\\
\\
and \\
\\
R C King\thanks{e-mail: R.C.King@maths.soton.ac.uk} \\
School of Mathematics, University of Southampton, \\
Southampton SO17 1BJ, England \\
}

\maketitle

\begin{abstract}
Alternating sign matrices with a U--turn boundary (UASMs) are a recent
generalization of ordinary alternating sign matrices.  Here we show that variations of
these matrices are in bijective correspondence with certain symplectic 
shifted tableaux that were recently introduced in the context of a symplectic
version of Tokuyama's deformation of Weyl's denominator formula. This 
bijection yields a formula for the weighted enumeration of UASMs. In
this connection use is made of the link between UASMs and certain square ice 
configuration matrices.
\end{abstract}

\vfill\eject

\section{Introduction}
\label{SectIntro}
Alternating sign matrices with a U--turn boundary (UASMs) first
appeared in a paper by Tsuchiya~\cite{T98} but
have been given a wider audience by Kuperberg~\cite{K02} and 
Propp~\cite{P02} (who called them half alternating sign matrices). 
In this paper we introduce a generalization of UASMs, called
$\mu$--UASMs, that combine the U--turn notion with the 
$\mu$--generalization of alternating sign matrices~(ASMs) due to 
Okada~\cite{O93}, where $\mu$ is a partition all of whose parts are
distinct.  We show that there exists a natural correspondence 
between $\mu$--UASMs and the symplectic shifted tableaux of shifted 
shape $\mu$ defined elsewhere~\cite{HK02}, and prove that this 
correspondence is a bijection.

There exists an important connection between ordinary alternating
sign matrices (ASMs) and square ice that was used to provide a second
proof of the alternating sign matrix conjecture by
Kuperberg~\cite{K96}. The square ice model involves two-dimensional
grids populated by frozen water molecules taking up any one of 
six configurations, see for example the work of Lieb~\cite{L67},
Bressoud~\cite{B99} and Lascoux~\cite{L99}. With a suitable choice of
boundary conditions this model can be linked to UASMs in a bijective 
manner~\cite{K02}. Here we extend this to the case of $\mu$--UASMs. 
To make this connection explicit it is convenient
to introduce square ice configuration matrices. These are then used
to provide both $x$ and $t$--weightings of $\mu$--UASMs, that are an 
exact counterpart to corresponding weightings of symplectic shifted 
tableaux. 

Thanks to the bijection between $\mu$--UASMs and symplectic 
shifted tableaux we are then able provide a general formula for the 
weighted enumeration of both $\mu$--UASMs and the UASMs
themselves. The latter correspond to the special cases 
$\mu=\delta=(n,n-1,\ldots,1)$ for some positive integer $n$. 
The most basic corollary of our result is 
\beq
 \sum_{UA\in {\cal UA}} 2^{\neg(UA)} = 2^{n^{2}},
\label{I.1}
\eeq
where ${\cal UA}$ is the set of $2n\times n$ UASMs and $\neg(UA)$ is 
the number of $-1$'s in $UA$. This result was conjectured by Propp~\cite{P02}
and proved by Eisenk\"olbl~\cite{E02} and independently by Chapman. It is also
derivable from Kuperberg~\cite{K02}.
More generally, we show that 
\beq
\sum_{UA\in {\cal UA}} t^{\ssi(UA)+\bar(UA)} (1+t)^{\neg(UA)} = (1+t)^{n^2},
\label{I.2}
\eeq
where $\ssi(UA)$ and $\bar(UA)$ are parameters, defined below, associated 
with each $UA\in{\cal UA}$.

The organisation of the paper is such that alternating sign matrices
and their U--turn manifestations are introduced in Section~2 and
symplectic shifted tableaux in Section~3. It is pointed out that the
latter may be viewed as being constructed from a sequence of ribbon
strips. It is this structure which is exploited in Section~4 to prove
the bijective nature of a map, $\Psi$, from each $sp(2n)$--standard shifted 
tableau $ST$ of shape specified by a partition $\mu$, all of whose
parts are distinct, to a $2n\times m$ $\mu$-UASM $UA=\Psi(ST)$, where
$m=\mu_1$,the largest part of $\mu$.

As a precursor to invoking two independent types of weighting of both
$\mu$--UASMs and symplectic shifted tableaux, square ice graphs and
the corresponding square ice configuration matrices are introduced in
Section~5. These configuration matrices then provide a natural way to
motivate and describe two different types of weighting of
$\mu$--UASMs, namely an $x$--weighting and a $t$--weighting. 
Corresponding weightings are then provided for symplectic shifted
tableaux. These latter weightings are those known to be relevant both 
to the character theory of $sp(2n)$~\cite{K76,KElS83} and to the
deformation of Weyl's denominator formula for $sp(2n)$~\cite{HK02}.
This is then exploited in Section~6 to provide a set of variously
weighted enumeration formulae for $\mu$--UASMs, symplectic shifted
tableaux and square ice configuration matrices. 

\section{Alternating Sign Matrices}
\label{SectASMs}
\addtocounter{equation}{-2}

Alternating sign matrices, ASMs, are square matrices all of whose
elements are $0$, $1$ or $-1$, such that the first and last non-zero 
entries of each row and column are $1$'s and the non-zero entries 
within each row and column alternate in sign. See, for example, 
the  $4\times 4$ ASM $A$ in equation (\ref{II.1}). Here and elsewhere we use $\ov1$ 
to denote $-1$.
\beq
A = \left[\begin{array}{rrrr}
0 & 1 & 0 & 0\\
1 & \ov1& 0 & 1\\
0 & 0 & 1 & 0\\
0 & 1 & 0 & 0\\
\end{array}\right].
\label{II.1}
\eeq
The number, $A(n)$, of $n\times n$ ASMs is described by the famous formula:
\beq
A(n)=\prod_{j=0}^{n-1} \frac{(3j+1)!}{(n+j)!}.
\label{II.2}
\eeq
The first proof of this formula was given by Zeilberger~\cite{Z96}. 
A second proof is due to Kuperberg~\cite{K96}, and a complete history 
is to be found in Bressoud~\cite{B99}.

Okada~\cite{O93} generalized ASMs by defining a set of $2n\times m$ 
$\mu$--alternating sign matrices, $\mu$--ASMs, associated with each 
partition $\mu=(\mu_1,\mu_2,\ldots,\mu_n)$
whose parts $\mu_j$ for $j=1,2,\ldots,n$ are all distinct and positive. 
These $\mu$--ASMs have properties similar to ordinary ASMs, but have 
column sums $1$ only in those columns indexed by $q=\mu_j$ for some 
$j$ and have column sums $0$ in all the other columns indexed by 
$q\neq\mu_j$ for any $j$. More formally, for each partition $\mu$
of length $\ell(\mu)=n$, all of whose parts are distinct, and for 
which $\mu_1\leq m$, an $n\times m$ matrix $A=(a_{iq})$ belongs
to the set ${\cal A}^\mu(n)$ of $n\times m$ $\mu$--ASMs if the
following conditions are satisfied:
\beq
\begin{array}{lll}
({\rm O1})\qquad &a_{iq}\in\{-1,0,1\} &\hbox{for $1\leq i\leq n$, $1\leq
q\leq m$};\cr
({\rm O2}) &\sum_{q=p}^m a_{iq}\in\{0,1\} 
&\hbox{for $1\leq i\leq n$, $1\leq p\leq m$}; \cr
({\rm O3}) &\sum_{i=j}^n a_{iq}\in\{0,1\} 
&\hbox{for $1\leq j\leq n$, $1\leq q\leq m$}; \cr 
({\rm O4}) &\sum_{q=1}^m a_{iq}=1 &\hbox{for $1\leq i\leq n$}; \cr 
({\rm O5}) &\sum_{i=1}^n a_{iq}= 
\bigg\{ \begin{array}{ll}
 1& \mbox{if } q=\mu_k \mbox{ for some $k$}\\ 
 0& \mbox{otherwise}
\end{array} 
&\hbox{for $1\leq q\leq m$, $1\leq k\leq n$}.\cr
\end{array}
\label{II.3} 
\eeq

The alternating sign matrices with a U--turn boundary, UASMs, are a 
variation on ordinary ASMs developed by Kuperberg~\cite{K02} after a 
paper of Tsuchiya~\cite{T98}. UASMs have an even number of rows. 
Each column of a UASM is of the same form as that of
an ordinary ASM. Each successive pair of rows of a UASM reading 
first from right to left across the top row of the pair and then 
from left to right across the bottom row of the pair is like a row 
of an ASM. Typically we have the $6\times 3$ UASM $UA$
\beq
UA = \left[\begin{array}{rrrr}
1 & 0 & 0\\
\ov1 & 1 & 0\\
0 & \ov1 & 1\\
1 & 0 & 0\\
0 & 1 & 0\\
0 & 0 & 0 \\
\end{array}\right].
\label{II.4} 
\eeq
The number, $A_{U}(2n)$, of UASMs of size $2n\times n$ is
\cite{K02}
\beq
A_U(2n)= 2^n (-3)^{n^2}
\prod_{ \genfrac{}{}{0pt}{}{1\leq i\leq 2n+1}{1\leq k\leq n} }
\frac{1+6k-3i}{2n+1+2k-i}. 
\label{II.5}
\eeq
Alternatively, thanks to their connection with
vertically symmetric ASMs (VSASMs) or flip symmetric ASMs (FSASMs), and a
recurrence relation for the number of the latter due to
Robbins~\cite{R00}, we have
\beq
A_U(2n)= A_U(2n-2) \biggl(\genfrac{}{}{0pt}{}{6n-2}{2n}\biggr)
\bigg/\biggl(\genfrac{}{}{0pt}{}{4n-2}{2n}\biggr).
\label{II.6}
\eeq
with $A_U(2)=2$. In either case we obtain:
\beq
\begin{array}{cccccccc}
n& 1&2&3&4&5&6&\cdots\cr A_U(2n)
&2&2^2\cdot3&2^3\cdot26&2^4\cdot646&2^5\cdot45885&2^6\cdot9304650&\cdots \cr
\end{array}
\label{II.7}
\eeq

Here we extend UASMs to the case of $\mu$--alternating sign matrices 
with a U--turn boundary. These were first defined in Hamel and 
King~\cite{HK02} in the context of deformations of Weyl's denominator 
formula for characters of the symplectic group and were called 
{\em sp($2n$)--generalised alternating sign matrices}.
\begin{Definition}
Let $\mu$ be a partition of length $\ell(\mu)=n$, all of whose parts 
are distinct, and for which $\mu_1\leq m$. Then the matrix 
$UA=(a_{iq})$ is said to belong to the set ${\cal UA}^\mu(2n)$ of 
$\mu$--alternating sign matrices with a U--turn boundary if it is
a $2n\times m$ matrix whose elements $a_{iq}$ satisfy the conditions: 
\beq
\begin{array}{lll}
({\rm UA1})\qquad &a_{iq}\in\{-1,0,1\} 
&\hbox{for $1\leq i\leq 2n$, $1\leq q\leq m$};\cr 
({\rm UA2})\qquad &\sum_{q=p}^m a_{iq}\in\{0,1\} 
&\hbox{for $1\leq i\leq 2n$, $1\leq p\leq m$}; \cr
({\rm UA3})\qquad &\sum_{i=j}^{2n} a_{iq}\in\{0,1\} 
&\hbox{for $1\leq j\leq 2n$, $1\leq q\leq m$}. \cr 
({\rm UA4})\qquad &\sum_{q=1}^m (a_{2i-1,q}+a_{2i,q})=1 
&\hbox{for $1\leq i\leq n$};\cr 
({\rm UA5})\qquad &\sum_{i=1}^{2n} a_{iq}=
  \bigg\{ \begin{array}{ll}
 1& \mbox{if } q=\mu_k \mbox{ for some $k$}\\ 
 0&\mbox{otherwise}
\end{array} 
 &\hbox{for $1\leq q\leq m$, $1\leq k\leq n$}.\cr
\end{array}
\label{II.8} 
\eeq
\label{Def2.1}
\end{Definition}

In the case for which $\mu=\delta=(n,n-1,\ldots,1)$ and $m=n$, for 
which (UA5) becomes $\sum_{i=1}^{2n}a_{iq}=1$ for $1\leq q\leq n$, 
this definition is such that the set ${\cal UA}^\delta(2n)$ coincides 
with the set of U--turn alternating sign matrices, UASMs, defined by 
Kuperberg~\cite{K02}. The more general case is exemplified for
the partition $\mu=(9,7,6,2,1)$ and $n=5$ by:
\beq
UA = \left[\begin{array}{rrrrrrrrr}
0 & 0 & 0 & 0 & 0 & 0 & 0 & \ov1 & 1 \\
1 & 0 & 0 & 0 & 0 & 0 & 0 & 0 & 0 \\
0 & 0 & 0 & 0 & 0 & 0 & 0 & 1 & 0 \\
\ov1 & 1 & 0 & \ov1 & 0 & 0 & 1 & 0 & 0 \\
1 & 0 & \ov1 & 1 & 0 & 0 & 0 & 0 & 0 \\
0 & 0 & 0 & \ov1 & 0 & 1 & 0 & 0 & 0 \\
0 & \ov1 & 0 & 1 & 0 & 0 & 0 & 0 & 0 \\
0 & 0 & 1 & 0 & 0 & 0 & 0 & 0 & 0 \\
\ov1 & 1 & 0 & 0 & 0 & 0 & 0 & 0 & 0 \\
1 & 0 & 0 & 0 & 0 & 0 & 0 & 0 & 0 \\
\end{array}\right]
\in {\cal UA}^\mu(2n)
\label{II.9}
\eeq
As can be seen the successive column sums reading from right to left are 
$1,1,0,0,0,1,1,0,1$, with the $1$'s appearing in columns $1,2,6,7,9$
specifying the parts of $\mu$. The individual row sums reading from
top to bottom are $0,1,1,0,1,0,0,1,0,1$ so that all the $U$--turn row 
sums for consecutive pairs of rows are $1$, as required.

In the proof of the bijection between $\mu$--UASMs and symplectic 
shifted tableaux in Section \ref{SectBiject} it will be useful to
refine the matrix $UA$. Any $\mu$--UASM $UA$ contains two types of 
zeros: zeros for which there is a nearest non-zero element to the 
right in the same row taking the value $1$ (positive zeros), and all
other zeros (negative zeros). We can then define a map $\phi$ from 
the matrix $UA$ to a signature matrix $\phi(UA)$, replacing
positive zeros and positive ones with plus signs, and negative zeros 
and negative ones with minus signs. It should be noted that there is 
no ambiguity in determining which zeros are positive and which are
negative, so that for each $\mu$--UASM $UA$ the signature matrix
$\phi(UA)$ is unique. Moreover, to recover $UA$ from $\phi(UA)$ by
means of the inverse map $\phi^{-1}$ it is only necessary in each row to
replace each right-most $+$ in a continuous sequence of $+$'s by $1$
and all others $+$'s by $0$, and the right-most $-$ of any continuous 
sequence of $-$'s by $-1$, provided that its immediate right-hand 
neighbour is $+$, and all other $-$'s by $0$.
This is illustrated in the case of our example~(\ref{II.9}) by
\beq
UA = \left[\begin{array}{rrrrrrrrr}
0 & 0 & 0 & 0 & 0 & 0 & 0 & \ov1 & 1 \\
1 & 0 & 0 & 0 & 0 & 0 & 0 & 0 & 0 \\
0 & 0 & 0 & 0 & 0 & 0 & 0 & 1 & 0 \\
\ov1 & 1 & 0 & \ov1 & 0 & 0 & 1 & 0 & 0 \\
1 & 0 & \ov1 &1 & 0 & 0 & 0 & 0 & 0 \\
0 & 0 & 0 & \ov1 & 0 & 1 & 0 & 0 & 0 \\
0 & \ov1 & 0 & 1 & 0 & 0 & 0 & 0 & 0 \\
0 & 0 & 1 & 0 & 0 & 0 & 0 & 0 & 0 \\
\ov1 & 1 & 0 & 0 & 0 & 0 & 0 & 0 & 0 \\
1 & 0 & 0 & 0 & 0 & 0 & 0 & 0 & 0 \\
\end{array}\right]
~\Rightarrow~\phi(UA) = \left[\begin{array}{rrrrrrrrr}
\m & \m & \m & \m & \m & \m & \m & \m & \p \\
\p & \m & \m & \m & \m & \m & \m & \m & \m \\
\p & \p & \p & \p & \p & \p & \p & \p & \m \\
\m & \p & \m & \m & \p & \p & \p & \m & \m \\
\p & \m & \m & \p & \m & \m & \m & \m & \m \\
\m & \m & \m & \m & \p & \p & \m & \m & \m \\
\m & \m & \p & \p & \m & \m & \m & \m & \m \\
\p & \p & \p & \m & \m & \m & \m & \m & \m \\
\m & \p & \m & \m & \m & \m & \m & \m & \m \\
\p & \m & \m & \m & \m & \m & \m & \m & \m \\
\end{array}\right]
\label{II.10}
\eeq

\section{Symplectic Shifted Tableaux}
\label{SecSymp}
\addtocounter{equation}{-10}

Symplectic shifted tableaux are variations on ordinary tableaux and 
were first introduced in~\cite{HK02} in the context of a symplectic 
version of Tokuyama's formula~\cite{T88} for the $t$--deformation of
Weyl's denominator formula.
A partition $\mu=(\mu_1,\mu_2,\ldots,\mu_n)$ is a weakly decreasing 
sequence of non-negative integers. The weight, $|\mu|$, of the
partition $\mu$ is the sum of its parts, and its length,
$\ell(\mu)\leq n$, is the number of its non-zero parts. Now
suppose all of the parts of $\mu$ are distinct. Define a shifted Young 
diagram $SF^\mu$ to be a set of $|\mu|$ boxes arranged in $\ell(\mu)$ 
rows of lengths $\mu_i$ that are left-adjusted to a diagonal line. 
More formally, 
$SF^{\mu}=\{(i,j)\,|\,1\leq i\leq\ell(\mu),i\leq j\leq \mu_i+i-1\}$.

For example, for $\mu=(9,7,6,2,1)$ we have
\beq
SF^{\mu}=\quad{\vcenter
{\offinterlineskip \halign{&\mystrut\vrule#&\mybox{\hss$\scriptstyle#$\hss}\cr
\hr{19}\cr 
& && && && && && && && && &\cr 
\hr{19}\cr 
\omit& && && && && && &&&& &\cr
\nr{2}&\hr{15}\cr 
\omit& &\omit& && && && &&&& && &\cr 
\nr{4}&\hr{13}\cr \omit& &\omit&
&\omit& && && &\cr 
\nr{6}&\hr{5}\cr \omit& &\omit& &\omit& &\omit& && &\cr
\nr{8}&\hr{3}\cr 
}}} 
\label{III.1}
\eeq
It should be noted that the parts of the partition
$\mu'=(\mu_1^\prime,\mu_2^\prime,\ldots,\mu_m^\prime)$, with
$m=\mu_1$, which is conjugate to $\mu$ specify the lengths 
of successive diagonals of $SF^\mu$. In the above example,
$\mu^\prime=(5,4,3,3,3,3,2,1,1)$. Quite generally, if all the
parts of $\mu$ are distinct, it follows that successive
parts of $\mu'$ differ by at most $1$. In fact, in such a case we have
\beq 
  \mu_{q+1}^\prime=
\left\{ \begin{array}{ll}
 \mu_q^\prime-1& \mbox{if } q=\mu_k \mbox{ for some $k$}\\ 
 \mu_q^\prime& \mbox {otherwise}
\end{array} \right. 
\label{III.2}
\eeq

Each symplectic shifted tableau, $ST$, is the result of filling
the boxes of $SF^\mu$ with integers from $1$ to $n$ and $\ov{1}$ 
to $\ov{n}$, ordered 
$\overline{1}<1 <\overline{2} <2 < \ldots <\overline{n} < n$, subject 
to a number of restrictions. We require a few more definitions. 
The {\em profile} of a shifted tableau is the sequence of entries 
on the main diagonal of the shifted tableau. Let $A$ be a
totally ordered set, or alphabet, and let $A^r$ be the set of 
all sequences $a=(a_1,a_2,\ldots,a_r)$ of elements of $A$ of length
$r$. Then the general set ${\cal ST}^\mu(A;a)$ is defined
to be the set of all standard shifted tableaux, $ST$, with respect 
to $A$, of profile $a$ and shape $\mu$, formed by placing an entry 
from $A$ in each of the boxes of $SF^\mu$ in such that the following 
five properties hold:
\beq
\begin{array}{cll}
({\rm S1})\quad& \eta_{ij} \in A 
& \quad\hbox{for all\ $(i,j)\in SF^\mu$};\\
({\rm S2})\quad& \eta_{ii} = a_i\in A 
& \quad\hbox{for all\ $(i,i)\in SF^\mu$};\\
({\rm S3})\quad& \eta_{ij} \leq \eta_{i,j+1}
& \quad\hbox{for all\ $(i,j), (i,j+1)\in SF^\mu$};\\
({\rm S4})\quad& \eta_{ij} \leq \eta_{i+1,j}
& \quad\hbox{for all\ $(i,j), (i+1,j)\in SF^\mu$};\\
({\rm S5})\quad& \eta_{ij} < \eta_{i+1,j+1} 
& \quad\hbox{for all\ $(i,j), (i+1,j+1)\in SF^\mu$}.
\end{array}
\label{III.3} 
\eeq

Informally, we may describe these tableaux as having shifted shape and 
as being filled with entries from $A$ with profile $a$ such that the 
entries are weakly increasing from left to right across each row and 
from top to bottom down each column, and strictly
increasing from top-left to bottom-right along each diagonal.

The set ${\cal ST}^\mu(n,\ov{n})$ of symplectic shifted tableaux is a 
specific instance of ${\cal ST}^\mu(A;a)$ given by:
\begin{Definition}
Let $\mu=(\mu_1,\mu_2,\ldots, \mu_n)$ be a partition of length
$\ell(\mu)=n$, all of whose parts are distinct, and let 
$A=[n,{\ov{n}}]=\{1,2,\ldots,n\}\cup\{\ov1,\ov2,\ldots,\ov{n}\}$ be 
subject to the order relations $\ov1<1<\ov2<2<\ldots<\ov{n}<n$. Then 
the set of all $sp(2n)$--standard shifted tableaux of shape $\mu$ is 
defined by:
\beq 
{\cal ST}^\mu(n,\ov{n})= \{S\in {\cal ST}^\mu(A;a)\,\vert\, 
A=[n,{\ov n}],~a\in[n,{\ov n}]^n ~\hbox{with $a_i\in\{i,{\ov i}\}$ 
for $i=1,2,\dots,n$} \},
\label{III.4} 
\eeq
where the entries $\eta_{ij}$ of each $sp(2n)$--standard shifted 
tableau $ST$ satisfy the conditions (S1)--(S5) of (\ref{III.3}). 
\label{Def3.1}
\end{Definition}

Continuing the above example with $n=5$ and $\mu=(9,7,6,2,1)$, we have
typically
\beq
ST=\ {\vcenter
{\offinterlineskip \halign{&\mystrut\vrule#&\mybox{\hss$\scriptstyle#$\hss}\cr
\hr{19}\cr
&\ov1&&1&&\ov2&&2&&\ov3&&\ov3&&\ov4&&4&&5&\cr 
\hr{19}\cr
\omit& &&\ov2&&\ov2&&2&&3&&\ov4&&\ov4&&4&\cr
\nr{2}&\hr{15}\cr
\omit& &\omit& &&3&&\ov4&&4&&4&&4&&4&\cr
\nr{4}&\hr{13}\cr
\omit& &\omit& &\omit& &&4&&4&\cr
\nr{6}&\hr{5}\cr
\omit& &\omit& &\omit& &\omit&&&\ov5&\cr \nr{8}&\hr{3}\cr
 }}}\ \in\ {\cal ST}^{97621}(5,\ov{5})
\label{III.5}
\eeq

Within each symplectic shifted tableau we can identify a further
construct, namely, a ribbon strip~\cite{HK02}.
\begin{Definition}
The ribbon strips $\str_k(ST)$ and $\str_{\ov{k}}(ST)$
consists of all boxes in the symplectic shifted tableau containing 
$k$ and ${\ov{k}}$, respectively, with no two such boxes on the 
same diagonal. Each ribbon strip may consist of one or more continuously
connected parts.
\label{Def3.2}
\end{Definition}

By way of example, for $ST$ as in (\ref{III.5}) $\str_4(ST)$ and 
$\str_{\ov4}(ST)$ take the form
\beq
\str_4(ST)=\ {\vcenter
{\offinterlineskip \halign{&\mystrut\vrule#&\mybox{\hss$\scriptstyle#$\hss}\cr
\nr{8}&\hr{3}\cr
\omit& &\omit& &\omit& &\omit& &&4&\cr
\nr{8}&\hr{3}\cr
\omit& &\omit& &\omit& &\omit& &&4&\cr
\nr{2}&\hr{9}\cr
\omit& &&4&&4&&4&&4&\cr
\hr{11}\cr
&4&&4&\cr
\hr{5}\cr
}}}
\qquad\qquad
\str_{\ov4}(ST)=\ {\vcenter
{\offinterlineskip \halign{&\mystrut\vrule#&\mybox{\hss$\scriptstyle#$\hss}\cr
\nr{6}&\hr{3}\cr
\omit& &\omit& &\omit& &&\ov4&\cr
\nr{4}&\hr{5}\cr
\omit& &\omit& &&\ov4&&\ov4&\cr
\hr{3}&\nr{1}&\hr{5}\cr
&\ov4&\cr
\hr{3}\cr
}}}.
\label{III.6}
\eeq

Each symplectic shifted tableaux is nothing other than a collection of 
ribbon strips nested or wrapped around one another so as to produce a diagram of 
standard shifted shape. 
It follows that each $ST\in{\cal ST}^\mu(n,\ov{n})$ may be encoded by
means of a map $\psi$ from $ST$ to a $2n\times m$ matrix $\psi(ST)$, with 
$m=\mu_1$, in which the rows of $\psi(ST)$, specified by $k$ and
$\ov{k}$ taken in reverse order from $n$ at the top to $\ov1$ at the
bottom, consist of a sequence of symbols $+$ or $-$ in the $q$th
column of $\psi(ST)$, counted from 1 on the left to $m$ on the right, 
indicating whether or not $\str_k(ST)$ and $\str_{\ov{k}}(ST)$, as 
appropriate, intersects the $q$th diagonal of $ST$, where diagonals 
are counted in the north-east direction starting from the main, first
diagonal to which the rows of $ST$ are left-adjusted.
Typically, applying $\psi$ to our example~(\ref{III.5}) for $ST$ gives
$\psi(ST)$ as shown:
\beq
ST=\ {\vcenter {\offinterlineskip
\halign{&\mystrut\vrule#&\mybox{\hss$\scriptstyle#$\hss}\cr
\hr{19}\cr
&\ov1&&1&&\ov2&&2&&\ov3&&\ov3&&\ov4&&4&&5&\cr 
\hr{19}\cr
\omit& &&\ov2&&\ov2&&2&&3&&\ov4&&\ov4&&4&\cr
\nr{2}&\hr{15}\cr
\omit& &\omit& &&3&&\ov4&&4&&4&&4&&4&\cr 
\nr{4}&\hr{13}\cr
\omit& &\omit& &\omit& &&4&&4&\cr
\nr{6}&\hr{5}\cr 
\omit& &\omit& &\omit& &\omit&&&\ov5&\cr 
\nr{8}&\hr{3}\cr
 }}}
~~\Rightarrow~~\psi(ST)= \left[\begin{array}{rrrrrrrrr}
\m & \m & \m & \m & \m & \m & \m & \m & \p \\
\p & \m & \m & \m & \m & \m & \m & \m & \m \\
\p & \p & \p & \p & \p & \p & \p & \p & \m \\
\m & \p & \m & \m & \p & \p & \p & \m & \m \\
\p & \m & \m & \p & \m & \m & \m & \m & \m \\
\m & \m & \m & \m & \p & \p & \m & \m & \m \\
\m & \m & \p & \p & \m & \m & \m & \m & \m \\
\p & \p & \p & \m & \m & \m & \m & \m & \m \\
\m & \p & \m & \m & \m & \m & \m & \m & \m \\
\p & \m & \m & \m & \m & \m & \m & \m & \m \\
\end{array}\right]~
\begin{array}{c} 5 \\ \ov5 \\  4 \\ \ov4 \\ 3 \\ \ov3 \\ 2 \\ \ov2 \\ 1 \\ \ov1 \\
\end{array}
\label{III.7}
\eeq

Clearly $\psi(ST)$ is uniquely determined by $ST$ and {\it vice versa}. 
The inverse map $\psi^{-1}$ from $\psi(ST)$ back to $ST$ is
accomplished by noting that the elements $+$ in
each column of $\psi(ST)$ simply signify by virtue of their row label,
$k$ or $\ov{k}$, 
those entries that appear in the corresponding diagonal of $ST$,
arranged in strictly increasing order.

The strips $\str_k(ST)$ and $\str_{\ov{k}}(ST)$, whose connected 
components are well represented by sequences of consecutive $+$'s in 
$\psi(ST)$, play a key role in establishing the bijection between
symplectic shifted tableaux and alternating sign matrices with a
U--turn boundary.

\section{The bijection}
\label{SectBiject}
\addtocounter{equation}{-7}

In Hamel and King \cite{HK02}, we derived a relationship between UASM
and symplectic shifted tableaux by first going through monotone
triangles. Here we prove the relationship directly. We will find it
useful to use the refinement of the UASM defined by $\phi$.

Since the image $\psi(ST)$ of $\psi$ acting on each symplectic shifted 
tableaux $ST$ is a matrix of $\pm$'s, the inverse $\phi^{-1}$ may
be applied to $\psi(ST)$ to give a matrix of $1$'s, $\ov1$'s and $0$'s,
which may or may not be a U--turn alternating sign matrix, $UA$. In 
fact the resulting matrix $\phi^{-1}\o\psi(ST)$ is always a U--turn 
alternating sign matrix, and it is shown below in Theorem~\ref{thmasm}
that the map $\Psi=\phi^{-1}\o\psi$ is a bijective mapping from 
${\cal ST}^\mu(n,\ov{n})$ to ${\cal UA}^\mu(2n)$.

In the case of our example, the outcome of this procedure mapping from $ST$ to
$\psi(ST)$, identifying $\psi(ST)$ with $\phi(UA)$, and then recovering 
$UA=\phi^{-1}\o\psi(ST)=\Psi(ST)$ is illustrated by: 
\beq 
{\vcenter {\offinterlineskip
\halign{&\mystrut\vrule#&\mybox{\hss$\scriptstyle#$\hss}\cr \hr{19}\cr
&\ov1&&1&&\ov2&&2&&\ov3&&\ov3&&\ov4&&4&&5&\cr
\hr{19}\cr
\omit& &&\ov2&&\ov2&&2&&3&&\ov4&&\ov4&&4&\cr
\nr{2}&\hr{15}\cr \omit& &\omit& &&3&&\ov4&&4&&4&&4&&4&\cr 
\nr{4}&\hr{13}\cr 
\omit& &\omit& &\omit& &&4&&4&\cr
\nr{6}&\hr{5}\cr
\omit& &\omit& &\omit& &\omit&&&\ov5&\cr 
\nr{8}&\hr{3}\cr }}}
\Longrightarrow \left[\begin{array}{rrrrrrrrr}
\m & \m & \m & \m & \m & \m & \m & \m & \p \\
\p & \m & \m & \m & \m & \m & \m & \m & \m \\
\p & \p & \p & \p & \p & \p & \p & \p & \m \\
\m & \p & \m & \m & \p & \p & \p & \m & \m \\
\p & \m & \m & \p & \m & \m & \m & \m & \m \\
\m & \m & \m & \m & \p & \p & \m & \m & \m \\
\m & \m & \p & \p & \m & \m & \m & \m & \m \\
\p & \p & \p & \m & \m & \m & \m & \m & \m \\
\m & \p & \m & \m & \m & \m & \m & \m & \m \\
\p & \m & \m & \m & \m & \m & \m & \m & \m \\
\end{array}\right]
\Longrightarrow \left[\begin{array}{rrrrrrrrr}
0 & 0 & 0 & 0 & 0 & 0 & 0 & \ov1 & 1 \\
1 & 0 & 0 & 0 & 0 & 0 & 0 & 0 & 0 \\
0 & 0 & 0 & 0 & 0 & 0 & 0 & 1 & 0 \\
\ov1 & 1 & 0 & \ov1 & 0 & 0 & 1 & 0 & 0 \\
1 & 0 & \ov1 & 1 & 0 & 0 & 0 & 0 & 0 \\
0 & 0 & 0 & \ov1 & 0 & 1 & 0 & 0 & 0 \\
0 & \ov1 & 0 & 1 & 0 & 0 & 0 & 0 & 0 \\
0 & 0 & 1 & 0 & 0 & 0 & 0 & 0 & 0 \\
\ov1 & 1 & 0 & 0 & 0 & 0 & 0 & 0 & 0 \\
1 & 0 & 0 & 0 & 0 & 0 & 0 & 0 & 0 \\
\end{array}\right]
\label{IV.1}
\eeq 
where the rows of the matrices are labelled from top to bottom
$n=5,\ov5,4,\ov4,3,\ov3,2,\ov2,1,\ov1$, and the columns from left 
to right $1,2,\ldots,9=m=\mu_1$.

\begin{Theorem}
Let $\mu=(\mu_1,\mu_2,\ldots,\mu_n)$ be a partition of length $\ell(\mu)=n$
whose parts are all distinct. Then 
the mapping $\Psi=\phi^{-1}\o\psi$ defines a bijection between the
set ${\cal ST}^\mu(n,\ov{n})$ of $sp(2n)$--standard shifted tableaux $ST$ of
shape $SF^\mu$, and the set ${\cal UA}^\mu(2n)$ of $2n\times m$
$\mu$--alternating sign matrices $UA$ with a U--turn boundary and
$m=\mu_1$.
\label{thmasm}
\end{Theorem}

\noindent{\bf Proof}: 
The Definition~\ref{Def3.1} of ${\cal ST}^\mu(n,\ov{n})$ ensures that
each $sp(2n)$--standard shifted tableau $ST$ satisfies the properties (S1)--(S5).
We need to show, in accordance with the Definition~\ref{Def2.1} of 
${\cal UA}^\mu(2n)$, that the properties (UA1)--(UA5)  hold for the matrix 
$UA=\Psi(ST)$ obtained from $ST$ by means of the map $\Psi$. 

First, it is obvious from the description of the mappings involved that 
the only possible matrix elements of $UA$ are $1$, $-1$, and $0$. 
Thus (UA1) holds.

Conditions (S3)--(S5) imply that each diagonal of $ST$ contains no
repeated entries, leading to the observation that $ST$ consists of a 
union of ribbon strips as described in Definition~\ref{Def3.2}. 
The map from $ST$ to the matrix $\psi(ST)$ is then such that reading
across each row of the matrix $\phi(UA)$ gives sequences of $+$'s 
corresponding to each connected component of the relevant ribbon
strip. The matrix $\psi(ST)$ is now to be identified with $\phi(UA)$ 
for some $UA$. The fact that the right-most $+$ of each sequence of 
consecutive $+$'s in $\phi(UA)$ is mapped to an element $1$ in $UA$, 
and that the right-most $-$ of each sequence of consecutive $-$'s is 
mapped to an element $-1$, provided that such a $-$ is followed by a
$+$, means that across each row of the resulting matrix $UA$ we have
non-zero entries $1$ and $-1$ that alternate in sign,
with the right-most non-zero entry always $1$. 
This implies the validity of (UA2).

To establish the U--turn nature of $UA$ it is necessary to invoke 
condition (S2) and the fact that $ST$ is standard only if the entry 
$\eta_{ii}=a_i$ in the $i$th box of the leading diagonal of $ST$ is 
either $i$ or $\ov{i}$. The map from $ST$ to $\psi(ST)$ is then
such that the elements in the first column of the $i$th and $\ov{i}$th 
rows are different, one is always $+$ and the other always $-$. 
Identifying $\psi(ST)$ with $\phi(UA)$, the first non-zero entries, if 
they exist, in the corresponding $i$th and $\ov{i}$th rows of $UA$
must also differ, one being $1$ and the other $-1$. This is
sufficient to show that the U--turn sequence obtained by reading 
across the $i$th row from right to left and then back along the 
$\ov{i}$th row from left to right is an alternating sequence of $1$'s 
and $-1$'s. The fact that in both rows the right-most non-zero element 
must be $1$ then ensures the validity of (UA4) since this U--turn 
alternating sign sequence begins and ends with $1$. If on the other 
hand either $i$ or $\ov{i}$ is not present in $ST$, then the 
corresponding row of $\psi(ST)$ will consist wholly of $-$'s, and 
identifying $\psi(ST)$ with $\phi(UA)$ leads to the conclusion that
the corresponding row of $UA$ consists solely of $0$'s, containing 
no non-zero elements and making no contribution to the U--turn
sequence. However, the other row of the pair $i$ and $\ov{i}$ in 
$\psi(ST)$ must start with a $+$ thereby ensuring that the first
non-zero entry in the corresponding row of $UA$ must be $1$. Since 
the last non-zero element is also $1$, the row sum is $1$ and the
U--turn condition (UA4) holds yet again.

To deal with (UA3), we consider the diagonals of $ST$. To this end the 
following schematic diagrams of various portions of the $q$th and 
$(q+1)$th diagonals of $ST$ will prove to be helpful.
\beq
D_1={\vcenter {\offinterlineskip
\halign{&\mystrut\vrule#&\mybox{\hss$#$\hss}\cr
\hr{5}\cr
&i&&b&\cr
\hr{7}\cr
\omit& &&a&&b&\cr
\nr{2}&\hr{7}\cr
\omit& &\omit& &&a&&b&\cr
\nr{4}&\hr{7}\cr
\omit& &\omit& &\omit& &&a&&b&\cr
\nr{6}&\hr{7}\cr
\omit& &\omit& &\omit& &\omit& &&a&&b&\cr
\nr{8}&\hr{5}\cr
\omit& &\omit& &\omit& &\omit& &\omit& &&j&\cr
\nr{10}&\hr{3}\cr}}}
~~~~~~~
D_2={\vcenter {\offinterlineskip
\halign{&\mystrut\vrule#&\mybox{\hss$#$\hss}\cr
\hr{3}\cr
&i&\cr
\hr{5}\cr
&a&&b&\cr
\hr{7}\cr
\omit& &&a&&b&\cr
\nr{2}&\hr{7}\cr
\omit& &\omit& &&a&&b&\cr
\nr{4}&\hr{7}\cr
\omit& &\omit& &\omit& &&a&&b&\cr
\nr{6}&\hr{7}\cr
\omit& &\omit& &\omit& &\omit& &&a&&j&\cr
\nr{8}&\hr{5}\cr}}}
 ~~~~~~~
D_3={\vcenter {\offinterlineskip
\halign{&\mystrut\vrule#&\mybox{\hss$#$\hss}\cr
\hr{9}\cr
\omit& &&a&&b&\cr
\nr{2}&\hr{7}\cr
\omit& &\omit& &&a&&b&\cr
\nr{4}&\hr{7}\cr
\omit& &\omit& &\omit& &&a&&b&\cr
\nr{6}&\hr{7}\cr
\omit& &\omit& &\omit& &\omit& &&a&&j&\cr
\nr{8}&\hr{5}\cr
}}}
\label{IV.2}
\eeq
In these diagrams the labels $i$ and $j$ are the actual entries 
in the corresponding boxes of $ST$, which may of course be barred or 
unbarred, while the rules (S3)--(S5) of (\ref{III.3}) are such that 
the actual entries of $ST$ in the boxes labeled by $a$ are all
distinct, as are those in the boxes labeled by $b$. Moreover, in
each case that we will consider each such entry $k$ will necessarily 
be such that $i<k<j$. We use the notation $n_a$ and $n_b$ to indicate 
the number of entries $a$ and $b$, respectively.

All elements $1$ in the $q$th column of the matrix $UA$ constructed 
from $ST$ by means of the map $\Psi$ correspond to connected
components of ribbon strips of $ST$ terminating in the $q$th diagonal, 
by virtue of their connection with right-most $+$'s in continuous 
sequences of $+$'s in the rows of $\psi(ST)=\phi(UA)$. Similarly all 
elements $-1$ in the $q$th column of the matrix $UA$ correspond to 
connected components of ribbon strips starting in the $q+1$th
diagonal, by virtue of their connection with the right-most $-$'s 
immediately preceding a $+$ in the rows of $\psi(ST)=\phi(UA)$. To see
that these non-zero elements in the $q$th column of $UA$ necessarily 
alternate in sign, consider two consecutive $1$'s and the
corresponding boxes on the $q$th diagonal of $ST$. In the schematic 
diagram $D_1$ above, these have been labeled by their entries $i$ and
$j$ (which could be barred or unbarred entries). They correspond
to the termination of connected components of the strips $\str_i(ST)$ 
and $\str_j(ST)$ in the $q$th diagonal of $ST$. All $n_a$ boxes on 
the $q$th diagonal between these $i$ and $j$ boxes, labeled in $D_1$ 
by $a$, must be labeled in $ST$ itself by $n_a$ distinct entries
$k$ with $i<k<j$. Similarly all $n_b$ boxes on the $(q+1)$th diagonal 
to the right of $i$ and above $j$, labeled in $D_1$ by $b$, must also 
be labeled in $ST$ by distinct entries $k$ with $i<k<j$. 
Since $n_b=n_a+1$ it follows that at least one $b$-label must be
distinct from all $a$-labels. If this label is $k$, then a connected 
component of $\str_k(ST)$ must start in the $(q+1)$th column with no 
component in the $q$th column. This leads in the $k$th row of
$\psi(ST)$ to a $-$ followed by a $+$, and hence to an
element $-1$ in the $q$th column of $UA$, between the two 
$1$'s associated with the boxes $i$ and $j$.

Similarly, between any two $-1$'s in the $q$th column of $UA$ there 
must exist an element $1$. The proof is based on the diagram $D_2$
above. The boxes labeled $i$ and $j$ in the $(q+1)$th diagonal of $ST$ 
specify the start of connected components of $\str_i(ST)$ and
$\str_j(ST)$ not present in the $q$th diagonal. Once again
in $D_2$ the set of $n_a$ boxes labeled by $a$ and the set of 
$n_b$ boxes labeled by $b$ must each have distinct labels $k$ with 
$i<k<j$ in $ST$ itself. Then $n_b=n_a-1$ so that there exists $k$, 
with $i<k<j$ such that $\str_k(ST)$ terminates in the $q$th column 
of $ST$, leading to a right-most $+$ in $\psi(ST)=\phi(UA)$ and hence 
to a $1$ in the $q$th column of $UA$ lying between the two $-1$'s 
associated with the $i$ and $j$ boxes of $ST$.

This is not sufficient to prove that (UA3) holds. It is necessary to 
prove further that the lowest non-zero entry in every column of $UA$
is $1$. The argument is very much as before. It should be noted that 
top and bottom are reversed in passing from $ST$ to $\psi(ST)$. We
consider the case of an entry $-1$ in the row corresponding to the 
label $j$ of the $q$th column of $UA$ and argue that there must exist 
an entry $1$ in the row corresponding to the label $k$ of the $q$th
column of $UA$ with $k<j$. This follows from the schematic diagram
$D_3$ which is truncated at its top end by the boundary of $ST$.
Again $n_b=n_a-1$ so that there exists $k$, with $k<j$ such that 
$\str_k(ST)$ terminates in the $q$th column of $ST$, leading to a 
right-most $+$ in $\psi(ST)=\phi(UA)$ and hence to a $1$ in the 
$q$th column of $UA$ lying below the $-1$ associated with the $j$ box of
$ST$. This applies to any element $-1$ in $UA$ so the lowest non-zero 
element of $UA$ must be $1$. In combination with the fact that, as we 
have proved, the signs of the non-zero elements are alternating in 
the columns of $UA$, this serves to complete the proof that (UA3) holds.

The final argument in respect of (UA5) is very similar. 
The relevant diagrams are as follows.
\beq
D_4={\vcenter {\offinterlineskip
\halign{&\mystrut\vrule#&\mybox{\hss$#$\hss}\cr
\hr{9}\cr
\omit& &&a&&b&\cr
\nr{2}&\hr{7}\cr
\omit& &\omit& &&a&&b&\cr
\nr{4}&\hr{7}\cr
\omit& &\omit& &\omit& &&a&&b&\cr
\nr{6}&\hr{7}\cr
\omit& &\omit& &\omit& &\omit& &&a&&b&\cr
\nr{8}&\hr{7}\cr
\omit& &\omit& &\omit& &\omit& &\omit& &&a&&\ast&\omit&\cr
\nr{10}&\hr{3}\cr
\omit& &\omit& &\omit& &\omit& &\omit& &\omit& &\cr
}}}
~~~~~~~
D_5={\vcenter {\offinterlineskip
\halign{&\mystrut\vrule#&\mybox{\hss$#$\hss}\cr
\hr{9}\cr
\omit& &&a&&b&\cr
\nr{2}&\hr{7}\cr
\omit& &\omit& &&a&&b&\cr
\nr{4}&\hr{7}\cr
\omit& &\omit& &\omit& &&a&&b&\cr
\nr{6}&\hr{7}\cr
\omit& &\omit& &\omit& &\omit& &&a&&b&\cr
\nr{8}&\hr{7}\cr
\omit& &\omit& &\omit& &\omit& &\omit& &&a&&b&&\phantom{b}&\omit&\cr
\nr{10}&\hr{7}\cr
\omit& &\omit& &\omit& &\omit& &\omit& &\omit& &&\ast&\omit&\cr
\nr{12}&\hr{1}\cr}}}
\label{IV.3}
\eeq

First we consider those diagonals $q$ of $ST$ which end with the 
right-most box of some row, that is those diagonals $q$ such that 
$q=\mu_k$ for some $k$ with $1\leq k\leq n$. As can be seen from 
the diagram $D_4$, in which the entry $\ast$ indicates an empty box 
just beyond the end of the $k$th row containing the final $a$ on the 
$q$th diagonal, the lengths $n_a=\mu_q^\prime$ and 
$n_b=\mu_{q+1}^\prime$ of the $q$th and $(q+1)$th
diagonals of $ST$, respectively, are such that $n_b=n_a-1$.
This is in accordance with (\ref{III.2}). 
Since the actual entries in $ST$ corresponding to the $n_a$ $a$'s are
all distinct, as are the entries corresponding to the $n_b$ $b$'s,
it follows that there exists precisely one more connected component
of the strips $\str_i(ST)$ that terminate in the $q$th diagonal of $ST$
than the number of connected components of strips $\str_j(ST)$ that 
start in the $(q+1)$th diagonal. Since it is the former that lead 
to all the $1$'s in the $q$th column of $UA$ and the latter to all 
the $-1$'s in the same column, the sum of the entries in this column 
must be $1$.

Similarly, we consider those diagonals $q$ of $ST$ which do not 
end with the right-most box of any row, that is those diagonals $q$ 
such that $q\neq\mu_k$ for any $k$. 
In this case the relevant schematic diagram is $D_5$ in which the
entry $\ast$ signifies an empty box just beyond the final $a$ in the $q$th
diagonal. Since this $a$ is not at the end of any row, it has a right
hand neighbour $b$, which must lie at the end of the $(q+1)$th
diagonal. As can be seen from $D_5$, in accordance with (\ref{III.2}), 
the lengths $n_a=\mu_q^\prime$ and $n_b=\mu_{q+1}^\prime$ of the $q$th 
and $(q+1)$th diagonals of $ST$, respectively, are such that $n_b=n_a$. 
This ensure that the number of connected components of strips
$\str_i(ST)$ that terminate in the $q$th diagonal of $ST$ is equal 
to the number of connected components of strips $\str_j(ST)$ that 
start in the $(q+1)$th diagonal. Once again, since it is the former 
that lead to all the $1$'s in the $q$th column of $UA$ and the latter
to all the $-1$'s in the same column, the sum of the entries in this column
must be $0$.

Taken together these last two results imply that (UA5) holds in all cases,
thereby completing the proof that for all 
$ST\in{\cal ST}^\mu(n,\ov{n})$ we have $UA=\Psi(ST)\in{\cal UA}^\mu(2n)$.

Reversing the argument, Definition~\ref{Def2.1} of ${\cal UA}^\mu(2n)$ 
ensures that each U--turn alternating sign matrix $UA$ satisfies the 
properties (UA1)--(UA5). 
We now need to show that these properties imply, in accordance with 
the Definition~\ref{Def3.1} of ${\cal ST}^\mu(n\ov{n})$, that 
$S=\Psi^{-1}(UA)$ satisfies (S1)--(S5), with 
$A=\{1,2,\ldots,n\}\cup\{\ov1,\ov2,\ldots,\ov{n}\}$ and
$a_i\in\{i,\ov{i}\}$ for $i=1,2,\ldots,n$.

First it should be noted that (UA1) guarantees the existence of 
$\phi(UA)=\psi(ST)$ as a matrix of $+$'s and $-$'s. The fact that 
$UA$ and hence $\psi(ST)$ is $2n\times m$, with rows labelled by the
elements of $A$, then ensures that (S1) holds,
since it is the row labels which determine the entries in $ST$.

The U--turn condition embodied in (UA2) and (UA4) then guarantees that
each pair of consecutive rows of $\phi(UA)$ counted from the bottom 
(or top) is such that one of the rows in the pair starts with a $+$
and the other with a $-$. In the case of the $i$th such pair, the row 
with $+$ in the first column of $\psi(ST)$ determines which one
of $i$ or $\ov{i}$ is the leading entry in the $i$th row of $ST$. 
This ensures that (S2) holds.

Thereafter, the fact that the entries of $ST$ are built up by adding
to the relevant diagonals all the $\ov1$'s, then all the $1$'s,
followed by all the $\ov2$'s, and so on, ensures that 
the ordering conditions (S3)--(S5) are 
automatically satisfied, provided that at every stage, after the
addition of all entries $\leq i$, the shape $SF^{\mu(i)}$
of the shifted sub-tableau, $S(i)$, obtained in this way is regular,
for all $i=\ov1,1,\ov2,\ldots,n$. By regular we mean that the lengths 
of the rows, left-adjusted as usual to the leading diagonal, are 
specified by means of a partition, $\mu(i)$, all of whose
non-vanishing parts are distinct.

To prove this we proceed by induction. First we consider the case
$i=\ov1$. The corresponding sub-tableau $S(i)=S(\ov1)$ is constructed
by adding to the empty diagram those boxes specified by the $+$'s
appearing in the bottom row of $\phi(UA)$ and filling them with
entries $\ov1$. Condition (UA3) with $j=2n$ ensures that the only
non-zero elements of the bottom row of $UA$ are $1$. The alternating
condition (UA2) with $i=2n$ then ensures that there is at most one
non-zero element $1$ in $UA$. If there is a $1$ in the $q$th column 
of $UA$, then the bottom row of $\phi(UA)=\psi(ST)$ consists of a
sequence of $q$ $+$'s followed by $(m-q)$ $-$'s. The procedure in
passing from $\psi(ST)$ to $ST$ then gives a sequence of $q$ boxes in its
top row, each containing the entry $\ov1$. Thus the shape of
$S(\ov1)$ is just $SF^q$, so that $\mu(\ov1)=(q)$, and
$SF^{\mu(\ov1)}$ is regular. If the bottom row of $UA$ contains no
entry $1$, then all its entries are $0$, the bottom row of
$\phi(UA)=\psi(ST)$ consists wholly of $-$'s and $S(\ov1)$ is empty.
By the U--turn condition the penultimate row of $UA$ will then contain at 
least one non-zero element. By the same argument as before using 
(UA2) and (UA3), there is only one such element and it must be a $1$. If
it lies in the $q$th column, then as before we find $S(1)$ consists
of a row of $q$ boxes with entries all equal to $1$. Once again we
have the shape $SF^q$, so that $\mu(1)=(q)$ and $SF^{\mu(1)}$ is
regular.

This serves to initialise the induction. Let $j$ be the element of
$A=[n,\ov{n}]$ that immediately precedes $i$ in the sequence
$\ov1,1,\ov2,2,\ldots,\ov{n},n$. We now assume that
$SF^{\mu(j)}$ corresponding to $S(j)$ is regular for some
$j\geq\ov1$ and consider the construction of $S(i)$ by adding to
$S(j)$ the strip $\str_i(ST)$, that is all the boxes with entries 
$i$ that are specified by the $+$'s appearing in the row of
$\phi(UA)$ labelled by $i$. As usual we concentrate on the 
$q$th and $(q+1)$th diagonals. By the induction hypothesis $S(j)$
is regular so that from (\ref{III.2}) the lengths of these diagonals, 
$n_a=\mu(j)_q^\prime$ and $n_b=\mu(j)_{q+1}^\prime$, are such that 
$n_b=n_a-1$ or $n_b=n_a$. 

We consider first the case $n_b=n_a-1$. The pair of elements in 
the $q$th and $(q+1)$th columns of $\phi(UA)$ can, at first sight,
be any one of the combinations $++$, $+-$, $-+$ and $--$. Taking 
these four possibilities in turn gives rise in $S(i)$ to $q$
and $(q+1)$th diagonals that are schematically of the form:
\beq
D_6={\vcenter {\offinterlineskip
\halign{&\mystrut\vrule#&\mybox{\hss$#$\hss}\cr
\hr{9}\cr
\omit& &&a&&b&\cr
\nr{2}&\hr{7}\cr
\omit& &\omit& &&a&&b&\cr
\nr{4}&\hr{7}\cr
\omit& &\omit& &\omit& &&a&&b&\cr
\nr{6}&\hr{7}\cr
\omit& &\omit& &\omit& &\omit& &&a&&b&\cr
\nr{8}&\hr{7}\cr
\omit& &\omit& &\omit& &\omit& &\omit& &&i&&i&\cr
\nr{10}&\hr{5}\cr
}}}
~~~
D_7={\vcenter {\offinterlineskip
\halign{&\mystrut\vrule#&\mybox{\hss$#$\hss}\cr
\hr{9}\cr
\omit& &&a&&b&\cr
\nr{2}&\hr{7}\cr
\omit& &\omit& &&a&&b&\cr
\nr{4}&\hr{7}\cr
\omit& &\omit& &\omit& &&a&&b&\cr
\nr{6}&\hr{7}\cr
\omit& &\omit& &\omit& &\omit& &&a&&b&\cr
\nr{8}&\hr{7}\cr
\omit& &\omit& &\omit& &\omit& &\omit& &&i&&\ast&\omit&\cr
\nr{10}&\hr{3}\cr
}}}
~~~
D_8={\vcenter {\offinterlineskip
\halign{&\mystrut\vrule#&\mybox{\hss$#$\hss}\cr
\hr{9}\cr
\omit& &&a&&b&\cr
\nr{2}&\hr{7}\cr
\omit& &\omit& &&a&&b&\cr
\nr{4}&\hr{7}\cr
\omit& &\omit& &\omit& &&a&&b&\cr
\nr{6}&\hr{7}\cr
\omit& &\omit& &\omit& &\omit& &&a&&b&\cr
\nr{8}&\hr{7}\cr
\omit& &\omit& &\omit& &\omit& &\omit& &\omit&\ast&&i&\cr
\nr{12}&\hr{3}\cr
}}}
~~~
D_9={\vcenter {\offinterlineskip
\halign{&\mystrut\vrule#&\mybox{\hss$#$\hss}\cr
\hr{9}\cr
\omit& &&a&&b&\cr
\nr{2}&\hr{7}\cr
\omit& &\omit& &&a&&b&\cr
\nr{4}&\hr{7}\cr
\omit& &\omit& &\omit& &&a&&b&\cr
\nr{6}&\hr{7}\cr
\omit& &\omit& &\omit& &\omit& &&a&&b&\cr
\nr{8}&\hr{7}\cr
\omit& &\omit& &\omit& &\omit& &\omit& &&\ast&\omit&\ast&\omit&\cr
%\nr{10}&\hr{1}\cr
}}}
\label{IV.4}
\eeq

As usual, $\ast$ signifies an empty box of $S(i)$, so that the diagrams 
$D_6$, $D_7$ and $D_9$ are regular in shape, but $D_8$ is irregular. 
It is therefore necessary to show that the shape $D_8$ never arises 
in passing from $UA$ to $S(i)$ by means of $\Psi^{-1}$. Since we have 
$n_b=n_a$ it follows, as in our previous discussion, that the $q$th 
column of $UA$ must contain the same number of $1$'s and $-1$'s below 
the row labeled by $i$. It then follows from (UA3) that any non-zero
element of $UA$ in the row labeled by $i$ and the $q$th column must be $1$. 
This excludes $D_8$ since it involves the start of a connected component 
of $\str_i(ST)$ in the $(q+1)$th diagonal. This can only arise from 
a pair $-+$ in $\phi(UA)$ and therefore an illegitimate element $-1$ in the row 
labelled by $i$ and the $q$th column of $UA$. Of course $D_7$ is allowed 
since it involves the end of a connected component of $\str_i(ST)$ in the
$q$th diagonal. This arises from a pair $+-$ in $\phi(UA)$ and therefore
a legitimate $1$ in the row labelled by $i$ and the $q$th column of $UA$.
Similarly $D_6$ and $D_9$ are both allowed since the corresponding pairs $++$
and $--$, respectively, of $\phi(UA)$ are associated with pairs of zeros of 
$UA$ itself, positive zeros in one case and negative zeros in the other.

The case of $n_b=n_a-1$ is similar. Once again the pair of elements in 
the $q$th and $(q+1)$th columns of $\phi(UA)$ can, at first sight,
be any one of the combinations $++$, $+-$, $-+$ and $--$. This time taking 
these four possibilities in turn gives rise in $S(i)$ to $q$ 
and $(q+1)$th diagonals that are schematically of the form:

\beq
D_{10}={\vcenter {\offinterlineskip
\halign{&\mystrut\vrule#&\mybox{\hss$#$\hss}\cr
\hr{9}\cr
\omit& &&a&&b&\cr
\nr{2}&\hr{7}\cr
\omit& &\omit& &&a&&b&\cr
\nr{4}&\hr{7}\cr
\omit& &\omit& &\omit& &&a&&b&\cr
\nr{6}&\hr{7}\cr
\omit& &\omit& &\omit& &\omit& &&a&&i&\cr
\nr{8}&\hr{5}\cr
\omit& &\omit& &\omit& &\omit& &\omit& &&i&\cr
\nr{10}&\hr{3}\cr
}}}
~~~
D_{11}={\vcenter {\offinterlineskip
\halign{&\mystrut\vrule#&\mybox{\hss$#$\hss}\cr
\hr{9}\cr
\omit& &&a&&b&\cr
\nr{2}&\hr{7}\cr
\omit& &\omit& &&a&&b&\cr
\nr{4}&\hr{7}\cr
\omit& &\omit& &\omit& &&a&&b&\cr
\nr{6}&\hr{5}\cr
\omit& &\omit& &\omit& &\omit& &&a&&\ast&\omit\cr
\nr{8}&\hr{5}\cr
\omit& &\omit& &\omit& &\omit& &\omit& &&i&\cr
\nr{10}&\hr{3}\cr
}}}
~~~
D_{12}={\vcenter {\offinterlineskip
\halign{&\mystrut\vrule#&\mybox{\hss$#$\hss}\cr
\hr{9}\cr
\omit& &&a&&b&\cr
\nr{2}&\hr{7}\cr
\omit& &\omit& &&a&&b&\cr
\nr{4}&\hr{7}\cr
\omit& &\omit& &\omit& &&a&&b&\cr
\nr{6}&\hr{7}\cr
\omit& &\omit& &\omit& &\omit& &&a&&i&\cr
\nr{8}&\hr{5}\cr
\omit& &\omit& &\omit& &\omit& &\omit& &&\ast&\omit&\cr
\nr{10}&\hr{1}\cr
}}}
~~~
D_{13}={\vcenter {\offinterlineskip
\halign{&\mystrut\vrule#&\mybox{\hss$#$\hss}\cr
\hr{9}\cr
\omit& &&a&&b&\cr
\nr{2}&\hr{7}\cr
\omit& &\omit& &&a&&b&\cr
\nr{4}&\hr{7}\cr
\omit& &\omit& &\omit& &&a&&b&\cr
\nr{6}&\hr{7}\cr
\omit& &\omit& &\omit& &\omit& &&a&&\ast&\omit&\cr
\nr{8}&\hr{3}\cr
\omit& &\omit& &\omit& &\omit& &\omit& &&\ast&\omit&\cr
\nr{10}&\hr{1}\cr
}}}
\label{IV.5}
\eeq 

Of these diagrams, $D_{10}$, $D_{12}$ and $D_{13}$ are regular and $D_{11}$
is irregular. In all these cases we have $n_b=n_a-1$ so that the $q$th column 
of $UA$ must contain one more $1$ than $-1$ below the row labeled 
by $i$. It follows from (UA3) that any non-zero element in the row labeled by 
$i$ must be $-1$. This excludes $D_{11}$ since this involves a connected component
of $\str_i(ST)$ ending in the $q$th diagonal of $S(i)$. This  
can only arise from a pair $+-$ in the row of $\psi(UA)$ specified by $i$, and
correspondingly from an illegitimate element $1$ in the $q$th column of 
this row of $UA$. Of course $D_{12}$ is allowed since this involves a 
connected component of $\str_i(ST)$ starting in the $(q+1)$th diagonal. This
necessarily arises from a pair $-+$ in the row of $\psi(UA)$ specified
by $i$, and correspondingly from a quite legitimate $-1$ in the $q$th
column of this row of $UA$. Similarly the cases $D_{10}$ and $D_{13}$ are allowed. 
They correspond to pairs of zeros in the $q$th and $(q+1)$th columns of 
the row of $UA$ specified by $i$, with one pair positive 
zeros and the other negative zeros.

Having excluded the irregular cases, represented by $D_{8}$ and $D_{11}$,
all the other cases are regular in that the length of the $(q+1)$th diagonal 
in $S(i)$ is either equal to the length of the $q$th diagonal or one less. 
This applies to all the diagonals of $S(i)$ so that its shape is given
by $SF^{\mu(i)}$ with $\mu(i)$ a partition, all of whose non-zero
parts are distinct. 

This completes the inductive argument passing from $S(j)$ of shape
$SF^{\mu(j)}$, with $\mu(j)$ a partition all of whose non-zero parts
are distinct, to $S(i)$ of shape $SF^{\mu(i)}$, with 
$\mu(i)$ also a partition all of whose non-zero parts are distinct. 
Since we have already dealt with the case $S(\ov{1})$ or $S(1)$, it follows 
that the shifted tableaux $S=S(n)$ is necessarily of shape $SF^{\mu(n)}$,
with $\mu(n)$ a partition all of whose non-zero parts are distinct.
The fact that $\mu(n)=\mu$ is then a consequence of (UA5). This can be
seen through a consideration of the diagram $D_4$, 
to be viewed now as part of $S(n))=\Psi^{-1}(UA)$. 
This diagram implies that for each $q$ such that 
$q=\mu(n)_i$ for some $i$, we have $n_b=n_a-1$. It follows that the number 
of $1$'s must have been one greater than the number of $-1$'s in the
$q$th column of $UA$. In accordance with the first case of (UA5), we
therefore have $q=\mu_j$ for some $j$, so that for each $i$ there
exists $j$ such that $\mu(n)_i=\mu_j$. Conversely, 
for each $q$ such that $q=\mu_j$ for some $j$ the number of $1$'s in
the $q$th column of $UA$ must be one greater than the number of $-1$'s in 
the same column. This implies that the $q$th and $(q+1)$th 
diagonals of $S(n)$ have the configuration given schematically in
$D_4$, so that $q=\mu(n)_i$ for some $i$. Thus for each $j$ there 
exists $i$ such that $\mu_j=\mu(n)_i$. 
Taken together these results ensure that $\mu(n)=\mu$, 
as required. Hence for all $UA\in{\cal UA}^\mu(2n)$ the conditions 
(S3)--(S5) apply to $S=\Psi^{-1}(UA)$. Having already established 
that (S1) and (S2) also apply, we can conclude that for all 
$UA\in{\cal UA}^\mu(2n)$ we have $S=\Psi^{-1}(UA)\in{\cal ST}^\mu(n,\ov{n})$.

This completes the proof of Theorem~\ref{thmasm} that $\Psi$ provides
a bijection between the $sp(2n)$--standard shifted tableaux 
$ST\in{\cal ST}^\mu(n,\ov{n})$ and 
the U--turn alternating sign matrices $UA\in{\cal UA}^\mu(2n)$ for all
partitions $\mu$ of length $\ell(\mu)=n$ whose parts are all distinct.

\section{Square ice}
\label{SecSquice}
\addtocounter{equation}{-5}

In order to exploit the above bijection to the full it is necessary to
add some $x$ and $t$--dependent weightings to both 
$UA\in{\cal UA}^\mu(2n)$ and $ST\in{\cal ST}^\mu(n,\ov{n})$. Although
some such weightings have already been provided~\cite{HK02}, rather
similar but not quite identical weightings may perhaps be best
motivated and described through the connection between $\mu$--UASMs, 
symplectic shifted tableaux and the square ice model that has proved 
to be such an invaluable tool in the study of alternating sign
matrices and their enumeration.

Square ice is a two dimensional grid that models the orientation of molecules
in frozen water, see for example Lieb~\cite{L67}, Bressoud~\cite{B99}, 
Lascoux~\cite{L99}.  In frozen water the model is such that each
individual molecule, consisting of two hydrogen atoms attached to an 
oxygen atom, takes up one of the $6$ possible orientations 
(the six vertex model) shown below.
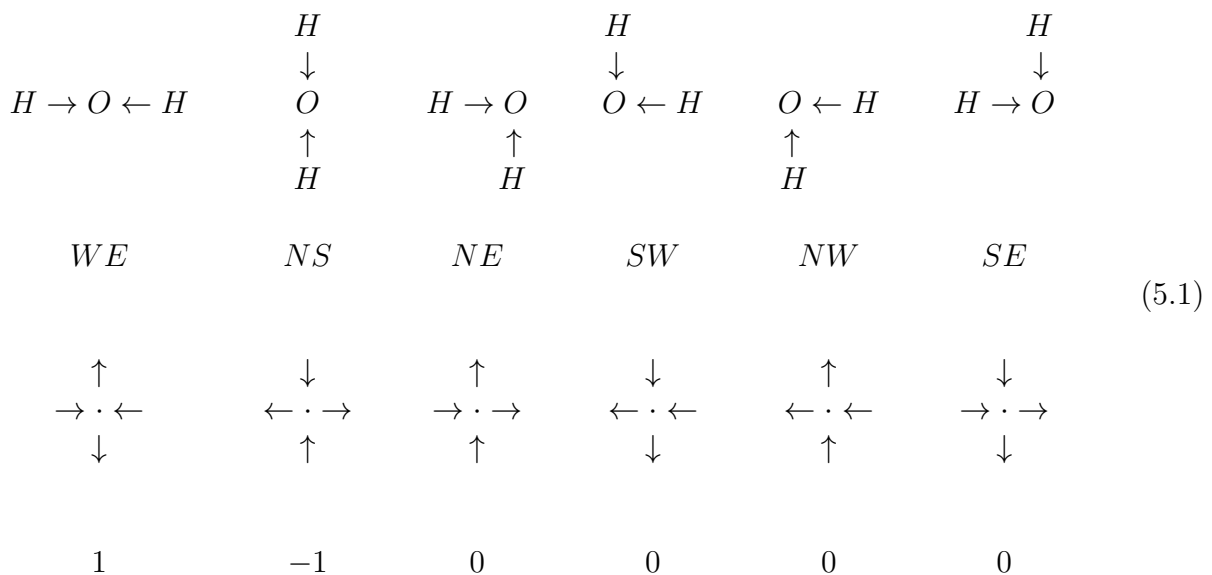
\beq
\begin{array}{ccccccccccc}
    &    &H   &    &        &    &H\,\phE&    &       &    &\phW\,H\cr
    &    &\da &    &        &    &\da\phE&    &       &    &\phW\da\cr
\HOH&\ \ &O   &\ \ &     \HO&\ \ &\OH    &\ \ &\OH    &\ \ &    \HO\cr
    &    &\ua &    &\phW \ua&    &       &    &\ua\phE&    &       \cr
    &    &H   &    &\phW \,H&    &       &    &H\,\phE&    &       \cr
\cr
 WE &    &NS  &    &NE      &    &SW     &    &NW     &    &SE     \cr
\cr
             &   &             &       &              &   &                &  &       &    &    \cr
    \ua      &   &\da          &       &    \ua         &  &\da             &  &   \ua    &    & \da\cr
\ra \cdot \la&\ \ & \la \cdot \ra&  \  \ &  \ra \cdot \ra& \ \ &\la \cdot \la&\ \ & \la \cdot \la&\ \ & \ra\cdot\ra\cr
    \da      &    &\ua           &       &  \ua          &    &   \da     &    &\ua&    &      \da \cr
             &    &              &       &              &    &            &    &      &    &       \cr
\cr
 1           &    &-1            &       &0            &     &0    &
    &0    &    &0    \cr
\end{array}
\label{V.1}
\eeq

As indicated in the second line of (\ref{V.1}),
the orientation of each molecule may be specified by 
giving the compass directions of the bonds linking each hydrogen atom 
to the oxygen atom. Thus WE represents a horizontal molecule, NS a 
vertical molecule and NE, SW, NW and SE molecules in which the
hydrogen bonds are mutually perpendicular. Alternatively,
each oxygen atom may be associated with a tetravalent vertex 
with two incoming and two outgoing edges as shown in the 
third line of (\ref{V.1}). At each vertex it is the incoming
edges that are associated with the hydrogen bonds displayed
in the first line of (\ref{V.1}).

%The map from the square ice 
%configuration, consisting of a grid of such molecules, to the square
%ice graph is one in which each oxygen atom is mapped to a tetravalent 
%vertex with two incoming and two outgoing edges as shown in 
%equation (\ref{V.1a}). At each vertex it is the incoming edges that 
%are associated with the hydrogen bonds.
%
%
%\beq
%\begin{array}{ccccccccccc}
%             &   &             &       &              &   &                &  &       &    &    \cr
%    \ua      &   &\da          &       &    \ua         &  &\da             &  &   \ua    &    & \da\cr
%\ra \cdot \la&\ \ & \la \cdot \ra&  \  \ &  \ra \cdot \ra& \ \ &\la \cdot \la&\ \ & \la \cdot \la&\ \ & \ra\cdot\ra\cr
%    \da      &    &\ua           &       &  \ua          &    &   \da     &    &\ua&    &      \da \cr
%            &    &              &       &              &    &            &    &      &    &       \cr
%\cr
% 1           &    &-1            &       &0            &     &0    &    &0    &    &0    \cr
%\end{array}
%\label{V.1a}
%\eeq

%\begin{figure}[htb]
%\begin{center}
%\includegraphics{sixvert1.eps}
%\caption{Six possible square ice configurations}
%\label{figsix}
%\end{center}
%\end{figure}

Square ice configurations~\cite{L67} consist of arrangements of the above 
molecules with an oxygen atom at each point of a square $n\times n$ grid
The corresponding square ice graph~\cite{L67} is one in which the
internal vertices sit at the grid points specified by the oxygen
atoms. The particular boundary conditions that correspond to ASM
were apparently first considered by Korepin~\cite{K82}.
As we have indicated all the internal vertices are tetravalent, 
with two incoming and two outgoing edges. The boundary vertices,
including corner vertices, are not usually drawn. Corner vertices have 
no edges. Non-corner boundary vertices are of valency one, but there 
may be boundary conditions on the edges linking them to the internal 
vertices. Conventionally, each left or right non-corner boundary
vertex has an edge pointing towards the adjacent internal vertex, 
while each top or bottom non-corner boundary vertex has an edge
pointing away from the adjacent internal vertex. 

Each such square ice configuration is then associated with an
alternating sign matrix. To construct the ASM one merely associates
each internal vertex of the type shown in equation (\ref{V.1}) with 
the corresponding matrix element $1$, $-1$ or $0$ indicated in the 
bottom line of (\ref{V.1}). The fact 
that the corresponding matrix is an ASM is a consequence of the
boundary conditions and the fact that each hydrogen atom is linked to 
just one oxygen atom. Using this association Kuperberg employed known 
results on square ice to provide a second proof of the alternating 
sign matrix conjecture \cite{K96}. 

This natural link between square ice and ordinary ASMs may be
generalized slightly so as to account for the U--turns and zero sum 
columns of our $\mu$--UASMs. It is only necessary to modify the
boundary conditions. A zero sum in column $q$ corresponds to a square
ice graph with incoming rather than outgoing edges at the top boundary 
in column $q$. A U--turn corresponds to either an outgoing left
boundary edge at row $2i-1$ and an incoming left boundary edge 
at row $2i$, or an incoming left boundary edge at row $2i-1$ and 
an outgoing left boundary edge at row $2i$ as shown in 
Figure~\ref{V.1b}. With these changes in boundary conditions 
we can map the six types of vertices to $1$'s, 
$-1$'s, and $0$'s exactly as before and produce a $\mu$--UASM $UA$. 
However, the $0$'s carry less information than is available in the 
square ice graph. At an intermediate stage in mapping from the 
square ice graph to a $\mu$--UASM it is helpful to map to a square 
ice configuration matrix, $CM$, whose matrix elements are just the 
labels WE, NS, NE, SW, NW, and SE attached to the six types of
vertex in (\ref{V.1}). To be precise, we adopt the following
\begin{Definition} 
Let $\mu=(\mu_1,\mu_2,\ldots,\mu_n)$ be a partition
of length $\ell(\mu)=n$, all of whose parts are distinct and with largest
part $\mu_1=m$. Then the configuration matrix $CM$ belongs to the set 
${\cal CM}^\mu(2n)$ if it is the image under the map of its vertices to
matrix elements in the set $\{WE,NS,NE,SW,NW,SE\}$ defined in (\ref{V.1})
of a square ice graph on a $2n\times m$ grid in which each internal 
vertex has two incoming and two outgoing edges, with all
right-hand edges incoming, all bottom edges outgoing, each left-hand 
pair of edges a U--turn with one edge incoming and one outgoing, 
and all top edges either outgoing or incoming according as the column 
number counted from the left is or is not equal to one of the parts 
of $\mu$. 
\label{Def5.1}
\end{Definition}

This is exemplified for the square ice graph of
Figure~\ref{V.1b} by the corresponding configuration matrix
$CM$ given in (\ref{V.2}). 

%\begin{figure}[htb]
%\begin{center}
%\includegraphics{squareice1.eps}
%\caption{A U--turn square ice graph on a $2n\times m$ grid with $n=5$, 
%$m=9$ and $\mu=(9,7,6,2,1)$. Only the internal vertices are shown.}
%\label{figsquareice}
%\end{center}
%\end{figure}

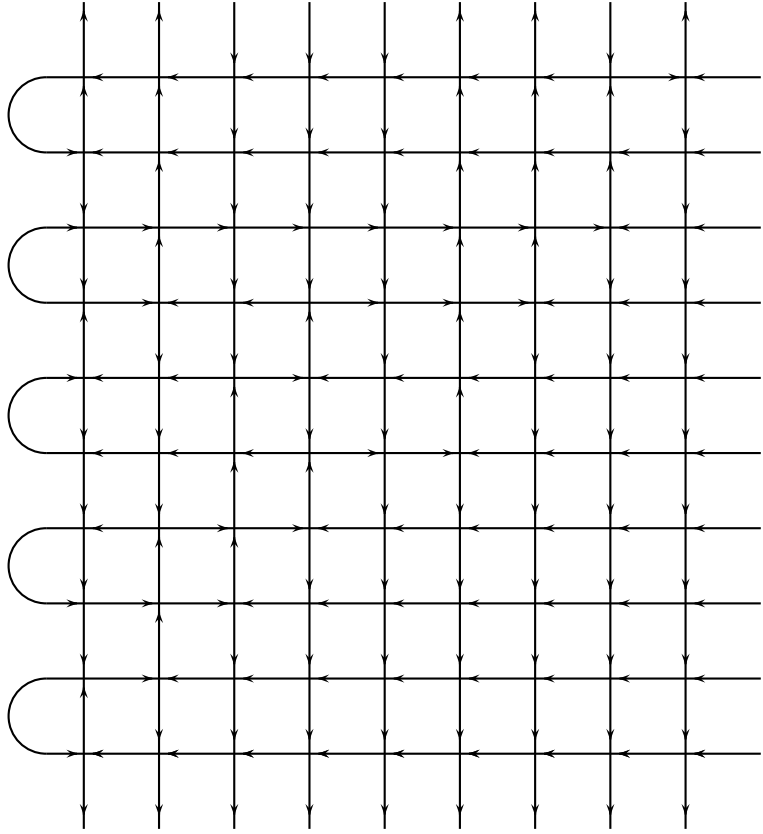
\begin{fullfigure}{V.1b}{Square ice with U--turn boundary}
\pspicture(-.5,-.5)(11,12)

\psline{->}(1.12,1)(1.11,1) \psline{->}(2.12,1)(2.11,1) \psline{->}(3.12,1)(3.11,1) 
\psline{->}(4.12,1)(4.11,1) \psline{->}(5.12,1)(5.11,1) \psline{->}(6.12,1)(6.11,1) 
\psline{->}(7.12,1)(7.11,1) \psline{->}(8.12,1)(8.11,1) \psline{->}(9.12,1)(9.11,1) 

\psline{->}(1.91,2)(1.92,2) \psline{->}(2.12,2)(2.11,2) \psline{->}(3.12,2)(3.11,2) 
\psline{->}(4.12,2)(4.11,2) \psline{->}(5.12,2)(5.11,2) \psline{->}(6.12,2)(6.11,2) 
\psline{->}(7.12,2)(7.11,2) \psline{->}(8.12,2)(8.11,2) \psline{->}(9.12,2)(9.11,2) 

\psline{->}(1.91,3)(1.92,3) \psline{->}(2.91,3)(2.92,3) \psline{->}(3.12,3)(3.11,3) 
\psline{->}(4.12,3)(4.11,3) \psline{->}(5.12,3)(5.11,3) \psline{->}(6.12,3)(6.11,3) 
\psline{->}(7.12,3)(7.11,3) \psline{->}(8.12,3)(8.11,3) \psline{->}(9.12,3)(9.11,3) 

\psline{->}(1.12,4)(1.11,4) \psline{->}(2.91,4)(2.92,4) \psline{->}(3.91,4)(3.92,4) 
\psline{->}(4.12,4)(4.11,4) \psline{->}(5.12,4)(5.11,4) \psline{->}(6.12,4)(6.11,4) 
\psline{->}(7.12,4)(7.11,4) \psline{->}(8.12,4)(8.11,4) \psline{->}(9.12,4)(9.11,4) 

\psline{->}(1.12,5)(1.11,5) \psline{->}(2.12,5)(2.11,5) \psline{->}(3.12,5)(3.11,5) 
\psline{->}(4.91,5)(4.92,5) \psline{->}(5.91,5)(5.92,5) \psline{->}(6.12,5)(6.11,5) 
\psline{->}(7.12,5)(7.11,5) \psline{->}(8.12,5)(8.11,5) \psline{->}(9.12,5)(9.11,5) 

\psline{->}(1.12,6)(1.11,6) \psline{->}(2.12,6)(2.11,6) \psline{->}(3.91,6)(3.92,6) 
\psline{->}(4.12,6)(4.11,6) \psline{->}(5.12,6)(5.11,6) \psline{->}(6.12,6)(6.11,6) 
\psline{->}(7.12,6)(7.11,6) \psline{->}(8.12,6)(8.11,6) \psline{->}(9.12,6)(9.11,6) 

\psline{->}(1.91,7)(1.92,7) \psline{->}(2.12,7)(2.11,7) \psline{->}(3.12,7)(3.11,7) 
\psline{->}(4.91,7)(4.92,7) \psline{->}(5.91,7)(5.92,7) \psline{->}(6.91,7)(6.92,7) 
\psline{->}(7.12,7)(7.11,7) \psline{->}(8.12,7)(8.11,7) \psline{->}(9.12,7)(9.11,7) 

\psline{->}(1.91,8)(1.92,8) \psline{->}(2.91,8)(2.92,8) \psline{->}(3.91,8)(3.92,8) 
\psline{->}(4.91,8)(4.92,8) \psline{->}(5.91,8)(5.92,8) \psline{->}(6.91,8)(6.92,8) 
\psline{->}(7.91,8)(7.92,8) \psline{->}(8.12,8)(8.11,8) \psline{->}(9.12,8)(9.11,8) 

\psline{->}(1.12,9)(1.11,9) \psline{->}(2.12,9)(2.11,9) \psline{->}(3.12,9)(3.11,9) 
\psline{->}(4.12,9)(4.11,9) \psline{->}(5.12,9)(5.11,9) \psline{->}(6.12,9)(6.11,9) 
\psline{->}(7.12,9)(7.11,9) \psline{->}(8.12,9)(8.11,9) \psline{->}(9.12,9)(9.11,9) 

\psline{->}(1.12,10)(1.11,10) \psline{->}(2.12,10)(2.11,10) \psline{->}(3.12,10)(3.11,10) 
\psline{->}(4.12,10)(4.11,10) \psline{->}(5.12,10)(5.11,10) \psline{->}(6.12,10)(6.11,10) 
\psline{->}(7.12,10)(7.11,10) \psline{->}(8.91,10)(8.92,10) \psline{->}(9.12,10)(9.11,10) 

\psline{->}(1.12,1)(1.11,1) \psline{->}(2.12,1)(2.11,1) \psline{->}(3.12,1)(3.11,1) 
\psline{->}(4.12,1)(4.11,1) \psline{->}(5.12,1)(5.11,1) \psline{->}(6.12,1)(6.11,1) 
\psline{->}(7.12,1)(7.11,1) \psline{->}(8.12,1)(8.11,1) \psline{->}(9.12,1)(9.11,1) 

\psline{->}(.91,1)(.92,1)
\psline{->}(.91,3)(.92,3)
\psline{->}(.91,6)(.92,6)
\psline{->}(.91,8)(.92,8)
\psline{->}(.91,9)(.92,9)

\psline{->}(1,.19)(1,.18) \psline{->}(1,1.88)(1,1.89)\psline{->}(1,2.19)(1,2.18)
\psline{->}(1,3.19)(1,3.18)\psline{->}(1,4.19)(1,4.18)\psline{->}(1,5.19)(1,5.18)
\psline{->}(1,6.88)(1,6.89)\psline{->}(1,7.19)(1,7.18)\psline{->}(1,8.19)(1,8.18)
\psline{->}(1,9.88)(1,9.89)\psline{->}(1,10.88)(1,10.89)

\psline{->}(2,.19)(2,.18) \psline{->}(2,1.19)(2,1.18)\psline{->}(2,2.88)(2,2.89)
\psline{->}(2,3.88)(2,3.89)\psline{->}(2,4.19)(2,4.18)\psline{->}(2,5.19)(2,5.18)
\psline{->}(2,6.19)(2,6.18)\psline{->}(2,7.88)(2,7.89)\psline{->}(2,8.88)(2,8.89)
\psline{->}(2,9.88)(2,9.89)\psline{->}(2,10.88)(2,10.89)

\psline{->}(3,.19)(3,.18) \psline{->}(3,1.19)(3,1.18)\psline{->}(3,2.19)(3,2.18)
\psline{->}(3,3.88)(3,3.89)\psline{->}(3,4.88)(3,4.89)\psline{->}(3,5.88)(3,5.89)
\psline{->}(3,6.19)(3,6.18)\psline{->}(3,7.19)(3,7.18)\psline{->}(3,8.19)(3,8.18)
\psline{->}(3,9.19)(3,9.18)\psline{->}(3,10.19)(3,10.18)

\psline{->}(4,.19)(4,.18) \psline{->}(4,1.19)(4,1.18)\psline{->}(4,2.19)(4,2.18)
\psline{->}(4,3.19)(4,3.18)\psline{->}(4,4.88)(4,4.89)\psline{->}(4,5.19)(4,5.18)
\psline{->}(4,6.88)(4,6.89)\psline{->}(4,7.19)(4,7.18)\psline{->}(4,8.19)(4,8.18)
\psline{->}(4,9.19)(4,9.18)\psline{->}(4,10.19)(4,10.18)

\psline{->}(5,.19)(5,.18) \psline{->}(5,1.19)(5,1.18)\psline{->}(5,2.19)(5,2.18)
\psline{->}(5,3.19)(5,3.18)\psline{->}(5,4.19)(5,4.18)\psline{->}(5,5.19)(5,5.18)
\psline{->}(5,6.19)(5,6.18)\psline{->}(5,7.19)(5,7.18)\psline{->}(5,8.19)(5,8.18)
\psline{->}(5,9.19)(5,9.18)\psline{->}(5,10.19)(5,10.18)

\psline{->}(6,.19)(6,.18) \psline{->}(6,1.19)(6,1.18)\psline{->}(6,2.19)(6,2.18)
\psline{->}(6,3.19)(6,3.18)\psline{->}(6,4.19)(6,4.18)\psline{->}(6,5.88)(6,5.89)
\psline{->}(6,6.88)(6,6.89)\psline{->}(6,7.88)(6,7.89)\psline{->}(6,8.88)(6,8.89)
\psline{->}(6,9.88)(6,9.89)\psline{->}(6,10.88)(6,10.89)

\psline{->}(7,.19)(7,.18) \psline{->}(7,1.19)(7,1.18)\psline{->}(7,2.19)(7,2.18)
\psline{->}(7,3.19)(7,3.18)\psline{->}(7,4.19)(7,4.18)\psline{->}(7,5.19)(7,5.18)
\psline{->}(7,6.19)(7,6.18)\psline{->}(7,7.88)(7,7.89)\psline{->}(7,8.88)(7,8.89)
\psline{->}(7,9.88)(7,9.89)\psline{->}(7,10.88)(7,10.89)

\psline{->}(8,.19)(8,.18) \psline{->}(8,1.19)(8,1.18)\psline{->}(8,2.19)(8,2.18)
\psline{->}(8,3.19)(8,3.18)\psline{->}(8,4.19)(8,4.18)\psline{->}(8,5.19)(8,5.18)
\psline{->}(8,6.19)(8,6.18)\psline{->}(8,7.19)(8,7.18)\psline{->}(8,8.88)(8,8.89)
\psline{->}(8,9.88)(8,9.89)\psline{->}(8,10.19)(8,10.18)

\psline{->}(9,.19)(9,.18) \psline{->}(9,1.19)(9,1.18)\psline{->}(9,2.19)(9,2.18)
\psline{->}(9,3.19)(9,3.18)\psline{->}(9,4.19)(9,4.18)\psline{->}(9,5.19)(9,5.18)
\psline{->}(9,6.19)(9,6.18)\psline{->}(9,7.19)(9,7.18)\psline{->}(9,8.19)(9,8.18)
\psline{->}(9,9.19)(9,9.18)\psline{->}(9,10.88)(9,10.89)
%\psline{->}(1,4.61)(1,4.62) \psline{->}(2,4.61)(2,4.62)
\psline(.5,1)(10,1) \psline(.5,2)(10,2)
\psline(.5,3)(10,3) \psline(.5,4)(10,4)
\psline(.5,5)(10,5) \psline(.5,6)(10,6)
\psline(.5,7)(10,7) \psline(.5,8)(10,8)
\psline(.5,9)(10,9) \psline(.5,10)(10,10)

\psline(1,0)(1,11) \psline(2,0)(2,11)
\psline(3,0)(3,11) \psline(4,0)(4,11)
\psline(5,0)(5,11) \psline(6,0)(6,11)
\psline(7,0)(7,11) \psline(8,0)(8,11)
\psline(9,0)(9,11)

\psarc(.5,1.5){.5}{90}{270}
\psarc(.5,3.5){.5}{90}{270}
\psarc(.5,5.5){.5}{90}{270}
\psarc(.5,7.5){.5}{90}{270}
\psarc(.5,9.5){.5}{90}{270}
\endpspicture
%\label{V.1b}
\end{fullfigure}

\beq
CM=\
\left[
\begin{array}{ccccccccc}
\ssc{NW}&\ssc{NW}&\ssc{SW}&\ssc{SW}&\ssc{SW}&\ssc{NW}&\ssc{NW}&\ssc{NS}&\ssc{WE}\\
\ssc{WE}&\ssc{NW}&\ssc{SW}&\ssc{SW}&\ssc{SW}&\ssc{NW}&\ssc{NW}&\ssc{NW}&\ssc{SW}\\
\ssc{SE}&\ssc{NE}&\ssc{SE}&\ssc{SE}&\ssc{SE}&\ssc{NE}&\ssc{NE}&\ssc{WE}&\ssc{SW}\\
\ssc{NS}&\ssc{WE}&\ssc{SW}&\ssc{NS}&\ssc{SE}&\ssc{NE}&\ssc{WE}&\ssc{SW}&\ssc{SW}\\
\ssc{WE}&\ssc{SW}&\ssc{NS}&\ssc{WE}&\ssc{SW}&\ssc{NW}&\ssc{SW}&\ssc{SW}&\ssc{SW}\\
\ssc{SW}&\ssc{SW}&\ssc{NW}&\ssc{NS}&\ssc{SE}&\ssc{WE}&\ssc{SW}&\ssc{SW}&\ssc{SW}\\
\ssc{SW}&\ssc{NS}&\ssc{NE}&\ssc{WE}&\ssc{SW}&\ssc{SW}&\ssc{SW}&\ssc{SW}&\ssc{SW}\\
\ssc{SE}&\ssc{NE}&\ssc{WE}&\ssc{SW}&\ssc{SW}&\ssc{SW}&\ssc{SW}&\ssc{SW}&\ssc{SW}\\
\ssc{NS}&\ssc{WE}&\ssc{SW}&\ssc{SW}&\ssc{SW}&\ssc{SW}&\ssc{SW}&\ssc{SW}&\ssc{SW}\\
\ssc{WE}&\ssc{SW}&\ssc{SW}&\ssc{SW}&\ssc{SW}&\ssc{SW}&\ssc{SW}&\ssc{SW}&\ssc{SW}\\
\end{array}
\right]
\label{V.2}
\eeq

The map $\chi$ from this configuration matrix $CM$ to the corresponding 
$\mu$--UASM $UA=\chi(CM)$ is then accomplished merely by setting 
WE and NS to $1$ and $-1$, respectively, and NE, SW, NW and SE all to
$0$. The example has been chosen so that the result is the matrix $UA$ 
appearing in (\ref{II.9}). It is not difficult to see that for all 
configuration matrices  $CM\in{\cal CM}^\mu(2n)$ we have 
$\chi(CM)\in{\cal UA}^\mu(2n)$. Moreover, $\chi$ is a bijection. 
The inverse map from $UA\in{\cal U}^\mu(2n)$ to 
$CM=\chi^{-1}(UA)\in{\cal CM}^\mu(2n)$ is such that the image under 
$\chi^{-1}$ of each matrix element $1$ and $-1$ of $UA$ is just WE and
NS, respectively. The images of the $0$'s are NE, SW, NW and SE 
according as their nearest non-zero neighbours to the right and 
below are $(1,1)$, $(\ov1,\ov1)$, $(\ov1,1)$ and $(1,\ov1)$,
respectively, where with some abuse of notation $\ov1$ is used to 
signify either $-1$ or the absence of any non-zero neighbour in 
the appropriate direction. These assignments are precisely what is 
required to ensure that there are no ambiguities in the directions 
of the edges at any vertex and that collectively they are consistent 
with the U--turn square ice conditions.

It is convenient to let $\we(CM)$, $\ns(CM)$, $\ne(CM)$, $\sw(CM)$, 
$\nw(CM)$, and $\se(CM)$ denote the total number of matrix elements 
of the configuration matrix $CM$ that are equal to WE, NS, NE, SW, NW 
and SE, respectively, and to refine this with subscripts $k$ and 
$\ov{k}$ if the count is restricted to the $(2n+1-2k)$th and 
$(2n+2-2k)$th rows, respectively. In addition we let $\ne_o(CM)$ 
and $\se_e(CM)$ denote the total number of matrix elements of $CM$
equal to NE in the odd rows counted from the top, and equal to SE 
in the even rows, and let $\wgt_e(CM)$ denote the total number
of matrix elements NE, SE and WE in the even rows. Thus
\beq
\begin{array}{rl}
\ne_o(CM)&=\sum_{k=1}^n\ \ne_k(CM);\cr
\se_e(CM)&=\sum_{k=1}^n\ \se_{\ov{k}}(CM);\cr
\wgt_e(CM)&=\sum_{k=1}^n \
(\ne_{\ov{k}}(CM)+\se_{\ov{k}}(CM)+\we_{\ov{k}}(CM)).\cr
\end{array}
\label{V.3'}
\eeq

The significance of these parameters and the fact that $\chi$ defines 
a bijection from $CM\in{\cal CM}^\mu(2n)$ to $UA\in{\cal UA}^\mu(2n)$ 
is that is that we may refer to the $0$'s of any such $UA=\chi(CM)$ 
as being NE, SW, NW or SE $0$'s if under $\chi^{-1}$ they 
map to NE, SW, NW or SE, respectively. Then, $\ne(CM)$, $\sw(CM)$, 
$\nw(CM)$ and $\se(CM)$ denote the numbers of such $0$'s in 
$UA=\chi(CM)$. In the same way the number of $1$'s and $-1$'s 
in $UA$ are given by $\we(CM)$ and $\ns(CM)$. Thus the configuration matrix 
$CM=\chi^{-1}(UA)$ is an alternative refinement of $UA$ to that
provided by the signature matrix $\phi(UA)$ exemplified in 
(\ref{II.10}). In fact the passage from $\phi(UA)$ to $CM$ is effected 
by replacing the right-most $+$ and right-most $-$ of any sequence of 
$+$'s and $-$'s in $\phi(UA)$ by $WE$ and $NS$, respectively, with the 
remaining $+$'s replaced by either $NE$ or $SE$ and the remaining
$-$'s by either $NW$ or $SW$ in accordance with the above rules 
regarding nearest non-zero neighbours of the corresponding $0$'s in
$UA$.

All this allows us to define various weightings and statistics on 
both $UA\in{\cal UA}^\mu(2n)$ and $ST\in{\cal ST}^\mu(n,\ov{n})$. 
First we define assign an $x$--weighting to each $\mu$--UASM. To 
this end let $m_k(UA)$ and $m_{\ov{k}}(UA)$ be the number of positive 
zeros and ones in the $k$th even and the $k$th odd row of $UA$, 
respectively, counted upwards from the bottom for $k=1,2,\ldots,n$. Then
\beq
 x^{\wgt(UA)}=x_1^{m_1(UA)-m_{\ov1}(UA)}\ x_2^{m_2(UA)-m_{\ov2}(UA)}
  \ \cdots\ x_n^{m_n(UA)-m_{\ov n}(UA)}.
\label{V.3}
\eeq
In our running example (\ref{II.10}) this gives
\beq
   x^{\wgt(UA)} = x_1^{1-1}\,x_2^{2-3}\,x_3^{2-2}\,x_4^{8-4}\,x_5^{1-1}
              =x_2^{-1}\, x_4^4.
\label{V.4}
\eeq
It should be noted that $m_k(UA)$ and $m_{\ov{k}}(UA)$ are just the number 
of $+$'s in the $(2n+1-2k)$th and $(2n+2-2k)$th rows of the signature
matrix $\phi(UA)$, respectively. 

Equivalently, in terms of the
configuration matrix $CM=\chi^{-1}(UA)$ we have
\beq
\begin{array}{rl}
  m_k(UA)&=m_k(CM)\quad\hbox{with}\quad 
m_k(CM)=\ne_k(CM)+\se_k(CM)+\we_k(CM);\\
  m_{\ov{k}}(UA)&=m_{\ov{k}}(CM)\quad\hbox{with}\quad 
m_{\ov{k}}(CM)=\ne_{\ov{k}}(CM)+\se_{\ov{k}}(CM)+\we_{\ov{k}}(CM),\\
\end{array}
\label{V.5}
\eeq
for $k=1,2,\ldots,n$. It follows that
\beq
    x^{\wgt(UA)}= x^{\wgt(CM)} \quad\hbox{with}\quad
x^{\wgt(CM)}=\prod_{k=1}^n\ x_k^{ m_k(CM)-m_{\ov{k}}(CM) }.
\label{V.6}
\eeq

There also exists a standard $x$--weighting of the
$sp(2n)$--symplectic shifted tableaux $ST$. To each entry $k$ or
$\ov{k}$ in $ST$ we associate a factor $x_k$ or $x_k^{-1}$. The
product of all these factors for $k=1,2,\ldots,n$ serves to define, 
as in \cite{HK02}, the $x$--weight of $ST$. Setting $m_k(ST)$ and 
$m_{\ov{k}}(ST)$ equal to the number of entries $k$ and $\ov{k}$, 
respectively, in $ST$ for $k=1,2\ldots,n$ we have
\beq
 x^{\wgt(ST)}=x_1^{m_1(ST)-m_{\ov1}(ST)}\ x_2^{m_2(ST)-m_{\ov2}(ST)}
  \ \cdots\ x_n^{m_n(ST)-m_{\ov n}(ST)}.
\label{V.7}
\eeq
In the example (\ref{III.5}) this gives
\beq
   x^{\wgt(ST)} = x_1^{1-1}\,x_2^{2-3}\,x_3^{2-2}\,x_4^{8-4}\,x_5^{1-1}
              =x_2^{-1}\, x_4^4.
\label{V.8}
\eeq

As can be seen from the bijective mapping from $ST$ to $UA=\Psi(ST)$ 
by way of $\psi(ST)=\phi(UA)$, illustrated in (\ref{II.10}), we have
\beq
\begin{array}{cl}
m_k(ST)&=m_k(UA) \quad\hbox{for $k=1,2,\ldots,n$;}\\ \\
m_{\ov{k}}(ST)&=m_{\ov{k}}(UA)\quad\hbox{for $k=1,2,\ldots,n$,}\\
\end{array}
\label{V.9}
\eeq
and hence
\beq
    x^{\wgt(ST)}=x^{\wgt(UA)}.
\label{V.10}
\eeq

In addition to the above $x$--weightings of both $UA\in{\cal UA}^\mu(2n)$
and $ST\in{\cal ST}^\mu(n,\ov{n})$, we can also assign $t$--weightings
to both $UA$ and $ST$. In dealing with $UA$ we require three statistics 
based on, but not quite identical to those introduced previously~\cite{HK02}. 
The first statistic, $\neg(UA)$, is defined to be the number of $-1$'s 
appearing in $UA$. The second statistic, $\bar(A)$, is defined to be
the total number of positive zeros and ones in the even rows of $UA$
counted from the top. This statistic can be read off most easily from 
$\phi(UA)$. For the third statistic we need the following:
\begin{Definition}
Let $UA$ be a $\mu$--UASM with matrix elements $a_{iq}$ for 
$1\leq i\leq 2n$ and $1\leq q\leq m$. Then $UA$ is said to have 
a {\em site of special interest}, an $\ssi$, at $(i,q)$ if:
\beq
\begin{array}{ll}
({\rm{SS1}})\qquad&\hbox{$a_{iq}=0$};\cr
({\rm{SS2}})\qquad&\hbox{$a_{ir}=1$ with $a_{ip}=0$ for $q<p<r\leq m$};\cr
({\rm{SS3}})\qquad&\hbox{either $i$ is odd and $a_{kq}=1$ with 
$a_{jq}=0$ for $i<j<k\leq 2n$},\cr
&\hbox{or $i$ is even and $a_{kq}=-1$ with $a_{jq}=0$ for $i<j<k\leq2n$},\cr
&\hbox{or $i$ is even and $a_{jq}=0$ for $i<j\leq 2n$}.\cr
\end{array}
\eeq
\label{Def5.2}
\end{Definition}

More graphically, each $\ssi$ is the site of a $0$ of $UA$ whose
nearest non-zero right hand neighbour is $1$, and whose nearest 
non-zero neighbour below the site is $1$ for a site in an odd row 
counted from the top and either $-1$ or non-existent for a site 
in an even row. With this definition, $\ssi(UA)$ is defined to be 
the number of sites of special interest in $UA$.

Once again it is perhaps easiest to read off these parameters $\neg(UA)$,
$\bar(UA)$ and $\ssi(UA)$ from the corresponding configuration matrix 
$CM=\chi^{-1}(UA)$. In terms of this matrix we have
\beq
\begin{array}{rl}
\neg(UA)&=\ns(CM)=\sum_{k=1}^n (\ns_k(CM)+\ns_{\ov{k}}(CM));\\ \\
\bar(UA)&=\sum_{k=1}^n\ (\ne_{\ov{k}}(CM)+\se_{\ov{k}}(CM)+ \we_{\ov{k}}(CM));\\ \\
\ssi(UA)&=\sum_{k=1}^n\ (\ne_k(CM)+\se_{\ov{k}}(CM)).\\
\end{array}
\label{V.11}
\eeq
 
In the example of (\ref{II.10}) we have $\neg(UA)=7$, $\bar(UA)=11$ 
and $\ssi(UA)=7$, where the seven sites of special interest are indicated 
by boldface ${\bf 0}$'s in $UA$, and by boldface $NE$'s and $SE$'s in 
$CM=\chi^{-1}(UA)$, as shown below in (\ref{V.12}).
\beq
UA=\left[\begin{array}{rrrrrrrrr}
0 & 0 & 0 & 0 & 0 & 0 & 0 & \ov1 & 1 \\
1 & 0 & 0 & 0 & 0 & 0 & 0 & 0 & 0 \\
0 & {\bf 0} & 0 & 0 & 0 & {\bf 0} &{\bf  0} & 1 & 0 \\
\ov1 & 1 & 0 & \ov1 & {\bf 0} & 0 & 1 & 0 & 0 \\
1 & 0 & \ov1 &1 & 0 & 0 & 0 & 0 & 0 \\
0 & 0 & 0 & \ov1 &{\bf 0 }& 1 & 0 & 0 & 0 \\
0 & \ov1 & {\bf 0} & 1 & 0 & 0 & 0 & 0 & 0 \\
{\bf 0 }& 0 & 1 & 0 & 0 & 0 & 0 & 0 & 0 \\
\ov1 & 1 & 0 & 0 & 0 & 0 & 0 & 0 & 0 \\
1 & 0 & 0 & 0 & 0 & 0 & 0 & 0 & 0 \\
\end{array}\right]
CM=\left[
\begin{array}{ccccccccc}
\ssc{NW}&\ssc{NW}&\ssc{SW}&\ssc{SW}&\ssc{SW}&\ssc{NW}&\ssc{NW}&\ssc{NS}&\ssc{WE}\\
\ssc{WE}&\ssc{NW}&\ssc{SW}&\ssc{SW}&\ssc{SW}&\ssc{NW}&\ssc{NW}&\ssc{NW}&\ssc{SW}\\
\ssc{SE}&\ssc{\bf{NE}}&\ssc{SE}&\ssc{SE}&\ssc{SE}&\ssc{\bf{NE}}&\ssc{\bf{NE}}
&\ssc{WE}&\ssc{SW}\\
\ssc{NS}&\ssc{WE}&\ssc{SW}&\ssc{NS}&\ssc{\bf{SE}}&\ssc{NE}&\ssc{WE}
&\ssc{SW}&\ssc{SW}\\
\ssc{WE}&\ssc{SW}&\ssc{NS}&\ssc{WE}&\ssc{SW}&\ssc{NW}&\ssc{SW}&\ssc{SW}&\ssc{SW}\\
\ssc{SW}&\ssc{SW}&\ssc{NW}&\ssc{NS}&\ssc{\bf{SE}}&\ssc{WE}&\ssc{SW}&\ssc{SW}&\ssc{SW}\\
\ssc{SW}&\ssc{NS}&\ssc{\bf{NE}}&\ssc{WE}&\ssc{SW}&\ssc{SW}&\ssc{SW}&\ssc{SW}&\ssc{SW}\\
\ssc{\bf{SE}}&\ssc{NE}&\ssc{WE}&\ssc{SW}&\ssc{SW}&\ssc{SW}&\ssc{SW}&\ssc{SW}&\ssc{SW}\\
\ssc{NS}&\ssc{WE}&\ssc{SW}&\ssc{SW}&\ssc{SW}&\ssc{SW}&\ssc{SW}&\ssc{SW}&\ssc{SW}\\
\ssc{WE}&\ssc{SW}&\ssc{SW}&\ssc{SW}&\ssc{SW}&\ssc{SW}&\ssc{SW}&\ssc{SW}&\ssc{SW}\\
\end{array}
\right]
\label{V.12}
\eeq

The $t$--weight to be attached to each element of $CM=\chi^{-1}(UA)$
can then be tabulated as follows
\beq
\begin{array}{ccc}
\quad\hbox{element}\quad&\quad\hbox{odd rows}\quad&\quad\hbox{even rows}\quad\\
WE&1&t\\
NS&1+t&1+t\\
NE&t&t\\
SW&1&1\\
NW&1&1\\
SE&1&t^2\\
\end{array}
\label{V.13}
\eeq
Comparison of (\ref{V.13}) with (\ref{V.11}) shows that this 
gives a total $t$--weight of
\beq
  t^{\ssi(UA)+\bar(UA)}\ (1+t)^{\neg(UA)}.
\label{V.14}
\eeq
Applying (\ref{V.13}) to (\ref{V.12}) gives the $t$--weighting
\beq
UA: (1+t)^7 \times \ 
\left[\begin{array}{l}
       1\ 1\ 1\ 1\ 1\ 1\ 1\ \ov1\ 1   \\
       t\ 1\ 1\ 1\ 1\ 1\ 1\ 1\ 1   \\
       1\ t\ 1\ 1\ 1\ t\ t\ 1\ 1   \\
       \ov1\ t\ 1\ \ov1\ t^2\ t\ t\ 1\ 1   \\
       1\ 1\ \ov1\ 1\ 1\ 1\ 1\ 1\ 1   \\
       1\ 1\ 1\ \ov1\ t^2\ t\ 1\ 1\ 1   \\
       1\ \ov1\ t\ 1\ 1\ 1\ 1\ 1\ 1   \\
       t^2\ t\ t\ 1\ 1\ 1\ 1\ 1\ 1   \\
       \ov1\ 1\ 1\ 1\ 1\ 1\ 1\ 1\ 1   \\
       t\ 1\ 1\ 1\ 1\ 1\ 1\ 1\ 1   \\
\end{array}\right] = t^{18}\ (1+t)^7.
\label{V.15}
\eeq

Turning now to the $t$--weighting of an $sp(2n)$--shifted tableau $ST$ 
it is convenient, in order to match contributions to the $t$--weight 
of $ST$ more precisely to the above contributions to the $t$--weight 
of $UA=\Psi(ST)$, to modify slightly our previous $t$--weighting of 
$sp(2n)$--standard shifted tableaux~\cite{HK02}. This is done as
follows. Each entry $k$ in $ST$ belongs to a ribbon strip $\str_k(ST)$ 
as in Definition \ref{Def3.2}. The $t$--weight of an entry $k$ is 
then defined to be $t$ if the entry immediately above this entry is 
also in $\str_k(ST)$, otherwise its $t$--weight is $1$. Similarly the 
$t$--weight of an entry $\ov{k}$ is defined to be $t^2$ if the entry 
immediately to its right is also in $\str_{\ov{k}}(ST)$, otherwise 
its $t$--weight is $t$. There is an additional $t$--weighting of
$(1+t)$ for every connected component of a strip $\str_k(ST)$ or 
$\str_{\ov{k}}(ST)$ that does not start on the main diagonal. 
In order to codify this, let $\str(ST)$ be the total number of
continuously connected components of all $\str_k(ST)$ and 
$\str_{\ov{k}}(ST)$ for $k=1,2,\ldots,n$, and let $\bar(ST)$ be the 
total number of barred entries in $ST$. In addition let
\beq
\var(ST)=\sum_{k=1}^n (\row_k(ST)-\con_k(ST)-\col_{\ov{k}}(ST)+\con_{\ov{k}}(ST)),
\label{V.16}
\eeq
where $\row_k(ST)$ and $\col_{\ov{k}}(ST)$ are the number of rows and
columns of $ST$ containing a $k$ and $\ov{k}$, respectively,
and $\con_k(ST)$ and $\con_{\ov{k}}(ST)$ are the number of
continuously connected components of $\str_k(ST)$ and 
$\str_{\ov{k}}(ST)$, respectively. This statistic $\var(ST)$
represents a measure of the upward steps in all 
$\str_k(ST)$ and the rightward steps in all $\str_{\ov{k}}(ST)$. 
In terms of the parameter $\hgt(ST)$ used in~\cite{HK02}, 
we have $\var(ST)=\hgt(ST)+\bar(ST)$.

For the strips of (\ref{III.6}) this $t$--weighting is illustrated by
\beq
\ {\vcenter
{\offinterlineskip \halign{&\mystrut\vrule#&\mybox{\hss$\scriptstyle#$\hss}\cr
\nr{10}&\hr{5}\cr 
\omit& &\omit& &\omit& &\omit& &\omit& &&1&&\cdot&\omit\cr
\nr{10}&\hr{3}\cr 
\omit& &\omit& &\omit& &\omit& &\omit& &&t&\cr
\nr{4}&\hr{9}\cr
\omit& &\omit& &&1&&1&&1&&t&\cr
\nr{2}&\hr{11}\cr
\omit& &&1&&t&\cr
\nr{2}&\hr{5}\cr
}}}
\qquad\qquad\qquad
{\vcenter
{\offinterlineskip \halign{&\mystrut\vrule#&\mybox{\hss$\scriptstyle#$\hss}\cr
\nr{8}&\hr{7}\cr
\omit& &\omit& &\omit& &\omit& &&t&&\cdot&\omit&\cdot&\omit\cr
\nr{6}&\hr{5}\cr
\omit& &\omit& &\omit&\cdot &&t^2&&t&\cr
\nr{2}&\hr{9}\cr
\omit& &&t&&\cdot\cr
\nr{2}&\hr{3}\cr
\omit& &&\cdot\cr
}}}\quad \times (1+t)^2.
\label{V.17}
\eeq
More generally, putting all such strips together we obtain the 
following $t$--weighting of $ST$ from (\ref{III.5}):
\beq
ST:\
{\vcenter
 {\offinterlineskip
 \halign{&\mystrut\vrule#&\mybox{\hss$\scriptstyle#$\hss}\cr
  \hr{19}\cr
  &t&&1&&t&&t&&t^2&\omit&t&&t&&1&&1&\cr
  \hr{5}&\nr{1}&\hr{1}&\nr{1}&\hr{5}&\nr{1}&\hr{1}&\nr{1}&\hr{3}\cr
  \omit& &&t^2&\omit&t&&1&&1&&t^2&\omit&t&&t&\cr
  \nr{2}&\hr{13}&\nr{1}&\hr{1}\cr
  \omit& &\omit& &&1&&t&&1&\omit&1&\omit&1&\omit&t&\cr
  \nr{4}&\hr{5}&\nr{1}&\hr{7}\cr
  \omit& &\omit& &\omit& &&1&\omit&t&\cr
  \nr{6}&\hr{5}\cr
  \omit& &\omit& &\omit& &\omit& &&t&\cr
  \nr{8}&\hr{3}\cr
 }}} \ \times \ (1+t)^7.
\label{V.18}
\eeq

As we have seen the bijection between $ST\in{\cal ST}^\mu(n,\ov{n})$ and
$UA\in{\cal UA}^\mu(2n)$ is such that the $(2n+1-2k)$th and
$(2n+2-2k)$th rows of $UA$ are determined by $\str_k(ST)$ and 
$\str_{\ov{k}}(ST)$, respectively, for $k=1,2,\ldots,n$. It is 
not difficult to see that each entry $k$ of weight $t$ corresponds 
to an NE $0$ of $UA$, while those of weight $1$ correspond either 
to a SE $0$ or to a WE entry $1$ if $k$ is the last entry of a 
connected component of $\str_k(ST)$. In the same way each entry 
$\ov{k}$ of weight $t^2$ corresponds to a SE $0$ of $UA$, while 
those of weight $t$ correspond either to a NE $0$ or to a WE entry 
$1$ if $\ov{k}$ is the last entry of a connected component of 
$\str_{\ov{k}}(ST)$. The additional weighting factors $(1+t)$ are 
associated with the NS $-1$'s of $UA$ since it is these $-1$'s that 
signal the start of a sequence of positive $0$'s ending in a $1$. 
In terms of the elements of the corresponding configuration matrix, 
$CM=\chi^{-1}(UA)$, arising from the square ice model we have
\beq
\begin{array}{rl}
\str(ST)&=\ns(CM)=\sum_{k=1}^n (\ns_k(CM)+\ns_{\ov{k}}(CM));\\ \\
\bar(ST)&=\sum_{k=1}^n\ (\ne_{\ov{k}}(CM)+\se_{\ov{k}}(CM)+\we_{\ov{k}}(CM));\\ \\
\var(ST)&=\sum_{k=1}^n\ (\ne_k(CM)+\se_{\ov{k}}(CM)).\\
\end{array}
\label{V.19}
\eeq
It follows from (\ref{V.11}) that for $UA=\Psi(ST)$ we have
\beq
\neg(UA)=\str(ST)-n, \quad \bar(UA)=\bar(ST), \quad \ssi(UA)=\var(ST).
\label{V.20}
\eeq

The coincidence of the $t$--weighting of $ST$ and $UA$ is exemplified 
in the case of our running example by
\beq
ST: (1+t)^7\ \times \
{\vcenter
 {\offinterlineskip
 \halign{&\mystrut\vrule#&\mybox{\hss$\scriptstyle#$\hss}\cr
  \hr{19}\cr
  &t&&1&&t&&t&&t^2&\omit&t&&t&&1&&1&\cr
  \hr{5}&\nr{1}&\hr{1}&\nr{1}&\hr{5}&\nr{1}&\hr{1}&\nr{1}&\hr{3}\cr
  \omit& &&t^2&\omit&t&&1&&1&&t^2&\omit&t&&t&\cr
  \nr{2}&\hr{13}&\nr{1}&\hr{1}\cr
  \omit& &\omit& &&1&&t&&1&\omit&1&\omit&1&\omit&t&\cr
  \nr{4}&\hr{5}&\nr{1}&\hr{7}\cr
  \omit& &\omit& &\omit& &&1&\omit&t&\cr
  \nr{6}&\hr{5}\cr
  \omit& &\omit& &\omit& &\omit& &&t&\cr
  \nr{8}&\hr{3}\cr
 }}}
\quad\Longleftrightarrow\quad
UA: (1+t)^7 \times \ \left[\begin{array}{l}
       0\ 0\ 0\ 0\ 0\ 0\ 0\ \ov1\ 1   \\
       t\ 0\ 0\ 0\ 0\ 0\ 0\ 0\ 0   \\
       1\ t\ 1\ 1\ 1\ t\ t\ 1\ 0   \\
       \ov1\ t\ 0\ \ov1\ t^2\ t\ t\ 0\ 0   \\
       1\ 0\ \ov1\ 1\ 0\ 0\ 0\ 0\ 0   \\
       0\ 0\ 0\ \ov1\ t^2\ t\ 0\ 0\ 0   \\
       0\ \ov1\ t\ 1\ 0\ 0\ 0\ 0\ 0   \\
       t^2\ t\ t\ 0\ 0\ 0\ 0\ 0\ 0   \\
       \ov1\ 1\ 0\ 0\ 0\ 0\ 0\ 0\ 0   \\
       t\ 0\ 0\ 0\ 0\ 0\ 0\ 0\ 0   \\
\end{array}\right]
\label{V.21}
\eeq
where, in particular, the 3rd and 4th rows of the $t$--weighting of
$UA$ are obtained from the $t$--weighting of $\str_4(ST)$ and
$\str_{\ov{4}}(ST)$ displayed in (\ref{V.18}). In contrast to
(\ref{V.15}) the NW and SW $0$'s of $UA$ have
been mapped to $0$ to indicate that they have no counterpart in
$ST$. In fact they correspond to the diagonals of $ST$ on which the
relevant strips have no box, as indicated for example by the $\cdot$'s
in (\ref{V.17}). 
Ignoring these $0$'s, the corresponding $t$--weight of both $ST$ and $UA$
is the product of all the displayed powers of $t$ together with the
seven factors $(1+t)$ arising from the seven continuously connected 
components of the strips of $ST$ that do not start on the main
diagonal, and equivalently from the seven $\ov1$'s of $UA$. It should
be noted that the sites of special interest in $UA$ correspond to the 
location of the $t$'s and $t^2$'s in the odd and even rows,
respectively, of the $t$--weighted version of $UA$.

It is perhaps worth summarising the $x$ and $t$--weighting by pointing 
out that in terms of the labelling used in $CM=\chi^{-1}(UA)$ the 
combined $x$ and $t$--weighting translates to
\beq
\begin{cases}
   (1+t)^{\ns(CM)}\ t^{\ne(CM)}\ (x_k)^{\ne(CM)+\se(CM)+\we(CM)}
   &\hbox{for row $2n+1-2k$};\cr
   \cr
   (1+t)^{\ns(CM)}\ t^{\se(CM)}\ (t\,x_k^{-1})^{\ne(CM)+\se(CM)+\we(CM)}
   &\hbox{for row $2n+2-2k$},\cr
\end{cases}
\label{V.22}
\eeq
for $k=1,2,\ldots,n$.

\section{Weighted Enumeration}
\label{SectWgt}
\addtocounter{equation}{-23}

Propp \cite{Pr00} has provided data for and made a number of
conjectures about the weighted enumeration of UASMs. 
Eisenk\"olbl~\cite{E02} has proved a number of these conjectures, and 
a number are derivable from Kuperberg~\cite{K02}. Here we delineate a 
new family of weighted enumerations of UASMs, and more generally of 
$\mu$--UASMs through their connection with symplectic shifted tableaux 
and show their overlap with the results of Propp.

To do this in the greatest generality, we shall also need the notion 
of ordinary symplectic tableaux \cite{K76,KElS83,S89}.
Let $\lambda=(\lambda_1,\lambda_2,\ldots,\lambda_r)$ be a partition 
of length $\ell(\lambda)=r\leq n$ and weight $|\lambda|$. Each such 
partition specifies a Young diagram $F^\lambda$ consisting of
$|\lambda|$ boxes arranged in $\ell(\lambda)$ rows of length
$\lambda_i$ that are left adjusted to a vertical line. For example, 
for $\lambda=(4,3,3)$ we have
\beq 
F^{\lambda}=\quad{\vcenter
{\offinterlineskip \halign{&\mystrut\vrule#&\mybox{\hss$\scriptstyle#$\hss}\cr
\hr{9}\cr 
& && && && &\cr 
\hr{9}\cr 
& && && &\cr
\hr{7}\cr 
& && && &\cr 
\hr{7}\cr 
}}} 
\label{VI.1} 
\eeq
Each symplectic tableau, $T$, of shape $\lambda$ is then the result of 
filling the boxes of $F^\lambda$ with integers from $1$ to $n$ and
$\ov{1}$ to $\ov{n}$, ordered
$\overline{1}<1<\overline{2}<2<\ldots<\overline{n}<n$, 
subject to a number of restrictions. This time let the profile of a
tableau be the sequence of entries obtained by reading down the first, 
left-most column.
Let $A$ be a totally ordered set, or alphabet, and let $A^r$ be the
set of all sequences $a=(a_1,a_2,\ldots,a_r)$ of elements of $A$ of
length $r$. Then the general set ${\cal T}^\lambda(A;a)$ is defined to 
be the set of all standard shifted tableaux, $ST$, with respect to
$A$, of profile $a$ and shape $\lambda$, formed by placing an entry 
from $A$ in each of the boxes of $F^\lambda$ in such that the
following four properties hold: 
\beq
\begin{array}{cll}
({\rm T1})\quad& \eta_{ij} \in A 
& \quad\hbox{for all\ $(i,j)\in F^\lambda$};\\
({\rm T2})\quad& \eta_{ii} = a_i\in A 
& \quad\hbox{for all\ $(i,1)\in F^\lambda$};\\
({\rm T3})\quad& \eta_{ij} \leq \eta_{i,j+1}
& \quad\hbox{for all\ $(i,j), (i,j+1)\in F^\lambda$};\\
({\rm T4})\quad& \eta_{ij} < \eta_{i+1,j}
& \quad\hbox{for all\ $(i,j), (i+1,j)\in F^\lambda$}.
\end{array}
\label{VI.2}
\eeq
These tableaux of shape $\lambda$ and profile $a$ have entries from $A$ that 
are weakly increasing from left to right across each row and are
strictly increasing from top to bottom down each column. 

The set ${\cal T}^\lambda(sp(2n))$ of all $sp(2n)$--standard tableaux
of shape $\lambda$ is a specific instance of ${\cal T}^\lambda(A;a)$
given by
\begin{Definition} Let
  $\lambda=(\lambda_1,\lambda_2,\ldots,\lambda_r)$ be 
a partition of length $\ell(\mu)=r\leq n$, and let 
$A=[n,{\ov n}] =\{1,2,\ldots,n\}\cup\{\ov1,\ov2,\ldots,\ov{n}\}$ 
be subject to the order relations $\ov1<1<\ov2<2<\ldots<\ov{n}<n$.
Then the set of all $sp(2n)$--standard  tableaux of shape $\lambda$ 
is defined by:
\beq
{\cal T}^\lambda(sp(2n))= 
     \{T\in {\cal T}^\lambda(A;a)\,\vert\, 
          A=[n,{\ov n}], a\in [n,{\ov n}]^r 
      \ \hbox{with $a_i\geq i$ for $i=1,2,\dots,r$}\},
\label{VI.3}
\eeq
where the entries $\eta_{ij}$ of each $sp(2n)$--standard tableau $T$ 
satisfy the conditions (T1)--(T4) of (\ref{VI.2}). 
\label{Def6.1}
\end{Definition}

Typically, for $n=5$ and $\lambda=(4,3,3)$ we have
\beq 
T=\quad{\vcenter
{\offinterlineskip \halign{&\mystrut\vrule#&\mybox{\hss$\scriptstyle#$\hss}\cr
\hr{9}\cr 
&\ov1&&\ov1&&1&&5&\cr 
\hr{9}\cr 
&2&&2&&4&\cr
\hr{7}\cr 
&4&&\ov5&&\ov5&\cr 
\hr{7}\cr 
}}} 
\quad\in\ {\cal T}^{433}(10).
\label{VI.4} 
\eeq

The symplectic Schur function \cite{K76,KElS83,S89}, 
which with a suitable interpretation of the indeterminates $x_i$ for 
$i=1,2,\ldots,n$ is the character of the irreducible representation of 
the Lie algebra $sp(2n)$ specified by $\lambda$, then takes the form 
\beq
 sp_\lambda(x) =  sp_\lambda(x_1,x_2,\ldots,x_n) =
      \sum_{T\in{\cal T}^\lambda(sp(2n))}\ x^{\wgt(T)},
\label{VI.5}
\eeq
where the sum is now over all $sp(2n)$--standard tableaux $T$ of shape
$\lambda$ and 
\beq
x^{\wgt(T)}= x_1^{m_1(T)-m_{\ov1}(T)}\ x_2^{m_2(T)-m_{\ov2}(T)}
  \ \cdots\ x_n^{m_n(T)-m_{\ov n}(T)},
\label{VI.6}
\eeq
with $m_k(T)$ and $m_{\ov k}(T)$ equal to the number of entries
$k$ and $\ov k$, respectively, in $T$. 

It is useful in the present context to generalise this by introducing 
some $t$-dependence and defining
\beq
sp_\lambda(x;t) =
      \sum_{T\in{\cal T}^\lambda(sp(2n))}\ t^{2\bar(T)}\ x^{\wgt(T)},
\label{VI.7}
\eeq
where $\bar(T)$ is the number of barred entries in $T$, that is 
\beq
  \bar(T)= \sum_{k=1}^n\ m_{\ov{k}}(T).
\label{VI.8}
\eeq
In our example (\ref{VI.4}) we have $\bar(T)=4$ and
\beq
   x^{\wgt(T)}=x_1^{1-2}\ x_2^{3-0}\ x_3^{0-0}\ x_4^{2-0}\ x_5^{0-2}
          =x_1^{-1}\ x_2^3\ x_4^2\ x_5^{-2}.
\label{VI.9}
\eeq

In this context the $t$--deformation of the denominator of Weyl's 
character formula for $sp(2n)$ takes the form
\beq
D_{sp(2n)}(x;t)=\!\!\prod_{1\leq i\leq n} x_i^{n-i+1}\!\!\!
   \prod_{1\leq i\leq n} (1+tx_i^{-2})\!\!\!\!
   \prod_{1\leq i<j\leq n} (1+tx_i^{-1}x_j)(1+tx_i^{-1}x_j^{-1}).
\label{VI.10}
\eeq

In our previous paper \cite{HK02} we derived the following extension 
to Tokuyama's formula \cite{T88} for the expansion of (\ref{VI.10}), 
namely
\begin{Theorem}
Let $\lambda$ be a partition into no more that $n$ parts and
let $\delta=(n,n-1,\ldots,1)$. Then
\begin{eqnarray}
D_{sp(2n)}(x;t)\ sp_\lambda(x;t)
     &=& \sum_{ST\in{\cal ST}^{\lambda+\delta}(n, \overline {n})}
     t^{\var(ST)+\bar(ST)}\ (1+t)^{\str(ST)-n}\ x^{\wgt(ST)},
\label{VI.11}
\end{eqnarray}
where the summation is taken over all $sp(2n)$--standard shifted 
tableaux $ST$ of shape $\mu=\lambda+\delta$.
\label{TheoST}
\end{Theorem}

Thanks to the bijection between $ST\in{\cal ST}^\mu(n,\ov{n})$ and 
$UA\in{\cal UA}^\mu(2n)$ and the equivalence between the $x$ and 
$t$--weightings of $ST$ and $UA$ this theorem can be recast in 
terms of $\mu$--UASMs as follows:
\begin{Theorem}
Let $\lambda$ be a partition into no more than $n$ parts,
let $\delta=(n,n-1,\ldots,1)$, and let $m=\lambda_1+n$. Then
\begin{eqnarray}
D_{sp(2n)}(x;t)\ sp_\lambda(x;t)
     &=& \!\!\! \sum_{UA\in{\cal UA}^{\lambda+\delta}(2n)}
    \!\!\! t^{\ssi(UA)+\bar(UA)}\ (1+t)^{\neg(UA)}\ x^{\wgt(UA)},
\label{VI.12}
\end{eqnarray}
where the summation is taken over all $2n\times m$ UASMs 
whose non-vanishing column sums are $1$ or $0$ according as the
column number is or is not a part of $\mu=\lambda+\delta$.
\label{TheoUA}
\end{Theorem}

Finally, in terms of the square ice configuration matrices we have
\begin{Theorem}\
Let $\lambda$ be a partition into no more than $n$ parts,
let $\delta=(n,n-1,\ldots,1)$ and let 
$m=\lambda_1+n$. Then
\begin{eqnarray}
D_{sp(2n)}(x;t)\ sp_\lambda(x;t)
     &=& \!\!\!\! \sum_{CM\in{\cal CM}^{\lambda+\delta}(2n)}
   \!\!\!\!\!  t^{\ne_o(CM)+\se_e(CM)+\wgt_e(CM)}\ (1+t)^{\ns(CM)}\ x^{\wgt(CM)},
\label{VI.13}
\end{eqnarray}
where the summation is taken over all $2n\times m$ U--turn 
square ice configuration matrices $CM$ whose top-most element in 
each column is either NW or NS if the column number is a part of 
$\mu=\lambda+\delta$ and is SW or WE otherwise.
\label{TheoCM}
\end{Theorem}

By setting $x_1=x_2=x_3=\ldots= x_n=1$ in these formulae we derive 
the following results
\begin{Corollary}\
Let $\lambda$ be a partition into no more than $n$ parts,
let $\delta$ be the partition $(n,n-1,\ldots,1)$, and let 
$m=\lambda_1+n$. Then
\begin{eqnarray}
(1+t)^{n^{2}}sp_\lambda(1;t)
&=&\sum_{ST\in{\cal ST}^{\lambda+\delta}(n, \overline {n})}
     t^{\var(ST)+\bar(ST)}\ (1+t)^{\str(ST)-n},\cr
&=&\!\!\!\! \sum_{UA\in{\cal UA}^{\lambda+\delta}(2n)}
    \!\!\! t^{\ssi(UA)+\bar(UA)}\ (1+t)^{\neg(UA)},\cr
&=&\!\!\!\! \sum_{CM\in{\cal CM}^{\lambda+\delta}(2n)}
   \!\!\!  t^{\ne_o(CM)+\se_e(CM)+\wgt_e(CM)}\ (1+t)^{\ns(CM)}.
\label{VI.14}
\end{eqnarray}
\label{Cor6.1}
\end{Corollary}

Specialising further to the case $t=1$ gives
\begin{Corollary}\
Let $\lambda$ be a partition into no more than $n$ parts,
let $\delta=(n,n-1,\ldots,1)$, and let 
$m=\lambda_1+n$. Then
\beq
   2^{n^{2}}sp_\lambda(1)=
\sum_{ST\in{\cal ST}^{\lambda+\delta}(n, \overline {n})}
     2^{\str(ST)-n}
=\!\!\! \sum_{UA\in{\cal UA}^{\lambda+\delta}(2n)}
    \!\!\! 2^{\neg(UA)}
=\!\!\! \sum_{CM\in{\cal CM}^{\lambda+\delta}(2n)}
   \!\!\!  2^{\ns(CM)}.
\label{VI.15}
\eeq
\label{Cor6.2}
\end{Corollary}

Here $sp_\lambda(1)=sp_\lambda(1;1)$ is the dimension of the 
irreducible representation of $sp(2n)$ specified by $\lambda$,
and it is known~\cite{W25,W26,ElSK79} that this is given by
\beq
sp_\lambda(1)=
 \prod_{1\leq i<j\leq n} \frac{\lambda_i -i-\lambda_j+j}{j-i} 
 \prod_{1\leq i \leq j\leq m}  
\frac{\lambda_i + \lambda_j + n -i -j +2}{n+2-i-j}.
\label{VI.16}
\eeq
However, as far as we know, no comparable product formula 
for $sp_\lambda(1;t)$ has yet been found.

Setting $\lambda=0$, so that $\mu=\delta=(n,n-1,\ldots,1)$ 
in Theorems \ref{TheoST}--\ref{TheoCM} gives
\begin{Theorem}
Let $\delta=(n,n-1,\ldots,1)$ and $m=n$. Then
\begin{eqnarray}
D_{sp(2n)}(x;t)
     &=& \sum_{ST\in{\cal ST}^{\delta}(n, \overline {n})}
     t^{\var(ST)+\bar(ST)}\ (1+t)^{\str(ST)-n}\ x^{\wgt(ST)}\cr
&=&\!\!\! \sum_{UA\in{\cal UA}^{\delta}(2n)}
    \!\!\! t^{\ssi(UA)+\bar(UA)}\ (1+t)^{\neg(UA)}\ x^{\wgt(UA)}\cr
&=&\!\!\! \sum_{CM\in{\cal CM}^{\delta}(2n)}
   \!\!\!  t^{\ne_o(CM)+\se_e(CM)+\wgt_e(CM)}\ (1+t)^{\ns(CM)}\ x^{\wgt(CM)},
\label{VI.17}
\end{eqnarray}
where the summations are taken over all $sp(2n)$--standard shifted 
tableaux $ST$ of shape $\delta$, all $2n\times n$ UASMs whose column 
sums are all $1$ and all $2n\times n$ U--turn square ice configuration 
matrices $CM$ whose top-most element in each column is either NW or NS. 
\label{TheoDelta}
\end{Theorem}

Finally, setting $x_1=x_2=\cdots=x_n=1$ in Corollaries \ref{Cor6.1} 
and \ref{Cor6.2} gives
\begin{Corollary}\
Let $\delta=(n,n-1,\ldots,1)$, then
\begin{eqnarray}
(1+t)^{n^{2}}&=&
\sum_{ST\in{\cal ST}^{\delta}(n, \overline {n})}
     t^{\var(ST)+\bar(ST)}\ (1+t)^{\str(ST)-n},\cr
&=&\!\!\! \sum_{UA\in{\cal UA}^{\delta}(2n)}
    \!\!\! t^{\ssi(UA)+\bar(UA)}\ (1+t)^{\neg(UA)},\cr
&=&\!\!\! \sum_{CM\in{\cal CM}^{\delta}(2n)}
   \!\!\!  t^{\ne_o(CM)+\se_e(CM)+\wgt_e(CM)}\ (1+t)^{\ns(CM)}.
\label{VI.18}
\end{eqnarray}
\label{Cor6.3}
\end{Corollary}
and 
\begin{Corollary}\
Let $\delta=(n,n-1,\ldots,1)$, then
\beq
   2^{n^{2}}=
\sum_{ST\in{\cal ST}^{\delta}(n, \overline {n})}
     2^{\str(ST)-n}
=\!\!\! \sum_{UA\in{\cal UA}^{\delta}(2n)}
    \!\!\! 2^{\neg(UA)}
=\!\!\! \sum_{CM\in{\cal CM}^{\delta}(2n)}
   \!\!\!  2^{\ns(CM)}.
\label{VI.19}
\eeq
\label{Cor6.4}
\end{Corollary}

\end{document}